\documentclass[hidelinks,onefignum,onetabnum]{siamart251104}
\usepackage{graphicx}
\usepackage{subfig}
\usepackage{color}
\usepackage{amsmath, amssymb, booktabs}
\usepackage{multirow}
\usepackage{xurl}
\usepackage{hyperref}

\usepackage{lipsum}
\usepackage{amsfonts}
\usepackage{graphicx}
\usepackage{epstopdf}
\usepackage{algorithmic}
\ifpdf
  \DeclareGraphicsExtensions{.eps,.pdf,.png,.jpg}
\else
  \DeclareGraphicsExtensions{.eps}
\fi

\newsiamremark{remark}{Remark}
\newsiamremark{hypothesis}{Hypothesis}
\crefname{hypothesis}{Hypothesis}{Hypotheses}
\newsiamthm{claim}{Claim}
\newsiamremark{fact}{Fact}
\crefname{fact}{Fact}{Facts}

\headers{on diagonalizable systems with random structure}{Yuan Zhang, Yutong Han, Yuanqing Xia, and Aming Li}

\title{An Example Article\thanks{Submitted to the editors DATE.
\funding{This work was funded by the Fog Research Institute under contract no.~FRI-454.}}}

\author{
Yuan Zhang\thanks{School of Automation, Beijing Institute of Technology, Beijing 100081, China 
  (\email{zhangyuan14@bit.edu.cn}, \email{hanyutong0105@163.com}).}
\and
Yutong Han\footnotemark[2]
\and
Yuanqing Xia\thanks{Zhongyuan University of Technology, Zhengzhou 450007, China, and  School of Automation, Beijing Institute of Technology, Beijing 100081, China (\email{xia\_yuanqing@bit.edu.cn}).}
\and
Aming Li\thanks{Center for Systems and Control, School of Advanced Manufacturing and Robotics, Peking University, Beijing 100871, China.
Research Center for Robotics, Peking University, Beijing 100871, China.
Center for Multi-Agent Research, Institute for Artificial Intelligence, Peking University, Beijing 100871, China
  (\email{amingli@pku.edu.cn}).}
}

\usepackage{amsopn}

  \ifpdf
\hypersetup{
  pdftitle={on diagonalizable systems with random structure},
  pdfauthor={Yuan Zhang, Yutong Han, Yuanqing Xia, and Aming Li}
}
\fi

\begin{document}
\title{on diagonalizable systems with random structure\thanks{Submitted to the editors DATE.
\funding{This work was supported in part by the National Key Research and Development Program of China under Grant 2025YFA1018800 and the National Natural Science Foundation of China under Grant 62373059. Corresponding authors: Yuanqing Xia and Aming Li.}}}
\maketitle
    \begin{abstract}Diagonalizability plays an important role in the analysis and design of multiple-variable systems. A structured matrix is called structurally diagonalizable if almost all of its numerical realizations—obtained by assigning real values to its free entries—are diagonalizable. Structural diagonalizability is useful in the verification and optimization of various structural system properties. In this paper, we study the asymptotic probability distribution of structural diagonalizability for structured systems whose system matrices are represented by directed Erdős–Rényi random graphs.  Leveraging a recently established graph-theoretic characterization of structural diagonalizability, we systematically analyze the distribution of structurally diagonalizable graphs under different edge-density regimes: for dense graphs, we prove that the system is almost always structurally diagonalizable; for graphs of medium density, we derive tight upper and lower bounds on the asymptotic probability of structural diagonalizability; and for extremely sparse graphs, we show that the asymptotic probability approaches $0$. The theoretical results are validated through extensive numerical simulations under varying numbers of vertices and connection probabilities.
    \end{abstract}
{\small{
\section{Introduction}
The analysis and design of complex systems often rely on simplifications of their dynamic models, where the diagonalizability of system matrices plays an important role.
Matrix diagonalizability, the property of being similar to a diagonal matrix, not only underpins controllability and stability analysis \cite{xue2021modal,chen1984linear}, modal decomposition \cite{sontag2013mathematical}, and controller design \cite{wang2003note}, but also directly determines whether the system dynamics can be decoupled, facilitating fault diagnosis and separation~\cite{chen2012robust}. By changing the coordinate system from the standard basis to an eigenbasis, diagonalizability can transform a multi-variable problem into a set of independent scalar problems.

In structured system theory, a structured matrix refers to a matrix whose entries are either fixed zeros or free values. A structured matrix (or a system represented by such a matrix) is defined as structurally diagonalizable if almost all its realizations, obtained by assigning real values to its free entries, are diagonalizable.
Recent research has demonstrated that structural diagonalizability can significantly simplify the verification of both structural output controllability and structural functional observability {\cite{montanari2022functional,montanari2023target,zhang2025generic,zhang2025structural}. Furthermore, for structurally diagonalizable systems, the associated actuator and sensor placement problems admit polynomial-time solutions; by contrast, these problems remain open for non-structurally diagonalizable systems ~\cite{zhang2023functional,zhang2025generic} (see Section \ref{sec-discussions} for details). Given that matrix diagonalizability is a generic property dominated by the zero-nonzero pattern of the corresponding matrices \cite{zhang2025generic}, the study of structural diagonalizability is significant not only in structured systems but also for numerical analysis.

Another control-theoretic application involving diagonalizability concerns spectral abscissa optimization for improving closed-loop stability \cite{BurkeLewisOverton2002}. Diagonalizability at a point implies local Lipschitzness of the spectral abscissa around that point, providing regularity of the corresponding optimization problem \cite{BurkeOverton1994}.   When the parameterized closed-loop matrix is compatible with
a structured-matrix model \cite{BurkeLewisOverton2002}, structural
diagonalizability indicates that almost every admissible realization is
diagonalizable and hence generically possesses this favorable local
regularity. 
Moreover, for large-scale networks comprised of linearly coupled subsystems, the (structural) diagonalizability of the underlying interconnection topology allows the controllability analysis of the entire network to be decomposed into two separate problems: the controllability of the global network topology and the inherent properties of the individual local subsystems \cite{zhang2021structural,xue2021modal}. This decomposition yields substantial simplifications for both the analysis and the design of such networked systems.}

However, topological structures of real-world systems—from sensor networks to neural circuits—are inherently uncertain due to stochastic link failures or probabilistic connections \cite{bullmore2009complex,yang2023reactivity}. This observation raises a fundamental research question: when a system's interconnection structure is itself random, to what extent are its dynamics inherently diagonalizable?

This paper investigates the probability that a structured matrix associated with a random directed graph is structurally diagonalizable. Existing research exhibits significant limitations in addressing this question. On one hand, structured systems theory for fixed interconnections \cite{lin1974structural, dion2003generic} has revealed the generic nature of system properties but has not addressed their probabilistic behavior under random structures. On the other hand, studies on the algebraic properties of random graphs have predominantly focused on spectral distributions \cite{wigner1958distribution} or specific attributes like structural controllability \cite{liu2011controllability}, while largely overlooking structural diagonalizability. 
This leaves a notable gap in the literature, as the structural diagonalizability of random structures has not yet been quantitatively characterized from a probabilistic perspective.

{\bf{Related work:}} Since Lin's seminal work on structural controllability \cite{lin1974structural}, the fundamental premise of this theoretical framework has been to infer dynamic properties exclusively from system zero-nonzero patterns. This idea was successfully extended to observability, stabilizability, and other system properties, culminating in a comprehensive theory of ``genericity", meaning these properties hold for almost all parameter realizations  \cite{reinschke1988multivariable,dion2003generic}. Recent research has shifted focus from the analysis of such generic properties toward their active design, robustness, and systematic optimization~\cite{Ramos2022AnOO}.

Research on random graph models initially focused on topological properties like connectivity and cliques \cite{erdos1960evolution, frieze2015introduction}. With the development of random matrix theory, research gradually delved into algebraic characterizations of graphs, such as the limiting distributions of adjacency matrix eigenvalues (e.g., the semicircle law) \cite{wigner1958distribution}. This perspective treats graphs as random matrices, whose spectral statistics become the focus \cite{tao2017random}. However, matrix diagonalizability represents a property distinct from spectral characteristics; it necessitates equality between geometric and algebraic multiplicities for each eigenvalue and demonstrates high sensitivity to specific numerical matrix entries. Although the spectra of random matrices are well-studied, the probability laws governing their diagonalizability remain challenging.

The convergence of these two research streams has catalyzed the probabilistic analysis of network structural properties. The seminal work by Liu et al. \cite{liu2011controllability} established the phase transition phenomenon for structural controllability in random graphs. They reveal that the number of driver vertices required to control the whole network is dominated by the degree distribution. Using random matrix theory, Rourke and Touri \cite{o2016conjecture} proved that the relative
number of controllable graphs compared to the total number of simple graphs on $n$ vertices approaches
one as $n$ tends to infinity. These studies offer deep insights into network controllability from a statistical perspective.

Notwithstanding these advances, existing research has not addressed probabilistic analysis of structural diagonalizability within the structured system framework. This fundamental property, governing a system's internal dynamic modes and directly impacting stability and modal decomposition, remains largely unexplored in the context of random networks.
To address this gap, we next introduce the random graph models used in this paper to formulate the probabilistic structural diagonalizability problem.

{\bf{Random graphs:}}
  The classic Erdős–Rényi (ER) random graph model is originally defined on undirected graphs \cite{erdos1960evolution}. To model networks with directional interactions, here we consider the directed graph extension version, further generalized by allowing self-loops to be parameterized independently. Next, we formally define the model $\mathcal{G}(n, p, q)$ and its special case $\mathcal{G}(n, p)$.

{\bf{The $\mathcal{G}(n, p, q)$ random graph model:}}
 For a graph with \( n \) labeled vertices, the connections are determined independently for each pair of vertices:

\begin{itemize}
    \item Non-loop directed edges: For any two distinct vertices \( u \) and \( v \) (where \( u \neq v \)), the probability that a directed edge exists from \( u \) to \( v \) is \( p \). Independently, the probability that a directed edge exists from \( v \) to \( u \) is also \( p \). The existence of each possible directed edge is an independent Bernoulli trial.
    \item Self-loop edges: For each vertex \( u \), the probability that a self-loop (a directed edge from \( u \) to itself) exists is \( q \). The existence of every self-loop is also independent.
    \item Independence: The existence of any edge (whether a non-loop directed edge or a self-loop) is independent of the existence of any other edge.
\end{itemize}

{\bf{The $\mathcal{G}(n, p)$ random graph model:}}
The model $\mathcal{G}(n, p, q)$ simplifies to the standard directed model \( \mathcal{G}(n, p) \) when the probability of self-loops equals the probability of non-loop edges, i.e., $q=p$. In this case, every possible directed edge appears independently with the same probability \( p \).}}

This paper presents the first systematic study of the probability of structural diagonalizability for random directed graphs. Based on the above-mentioned random graph models \(\mathcal{G}(n,p,q)\) and \(\mathcal{G}(n,p)\), we aim to characterize the probability of structural diagonalizability and uncover its evolution with edge density \(p\) and $q$. To this end, we first establish some easily verifiable necessary or sufficient graph-theoretic conditions for structural diagonalizability.  Leveraging these characterizations,  we derive analytical bounds and asymptotic properties for the probability of structural diagonalizability under random graph models \(\mathcal{G}(n,p,q)\) and \(\mathcal{G}(n,p)\). Specifically, we reveal that for dense graphs, the system is almost always structurally diagonalizable; for graphs of medium density, we establish tight upper and lower bounds on the asymptotic probability of structural diagonalizability; and for extremely sparse graphs, we show that this asymptotic probability approaches~$0$.

\section{Preliminaries: Structurally diagonalizable directed graphs}
\label{sec:main}
This section presents basic concepts in graph theory and establishes the framework for analyzing structural diagonalizability in linear structured systems.\par

Let $ n \in \mathbb{Z}^+ $ be a positive integer. For each pair $(i, j)$, the matrix $ e_{ij} $ is defined as an $ n \times n $ matrix where the entry in the $ i $-th row and $ j $-th column is a free variable (i.e., an arbitrary real number), while all other entries are zero. Let $ E \subseteq \{1, 2, \dots, n\} \times \{1, 2, \dots, n\} $. The \emph{zero pattern} $ Z_E $ is the linear subspace of $ \mathbb{R}^{n \times n} $ spanned by the basis matrices $ \{ e_{ij} : (i,j) \in E \} $,
\begin{equation}
    Z_{E}=\{e_{ij}\mid (i,j)\in E\}.
\end{equation}

A zero pattern defines a system structure encoded by a graph with edge set $ E $ (cf. \cite{belabbas2022stable}). It describes the system's connectivity through a {\emph{structured matrix}} $ \bar{A}\in \{0,*\}^{n\times n}$, defined as follows:
\begin{itemize}
    \item $ \bar{A}_{ij} = 0 $: Fixed zero (no interaction from vertex $ i $ to $ j $).
    \item $ \bar{A}_{ij} = \ast $: Free parameter (interaction exists, i.e., $ (i,j) \in E $).
\end{itemize}A realization of $\bar A$ is a matrix obtained by assigning real numbers to the free parameters of~$\bar A$ while preserving all fixed zero entries.

\begin{figure}[htbp]
\centering
\subfloat[$G_{\bar{A}_{\alpha}}$]{\label{fig:0a}\includegraphics[width=0.32\textwidth]{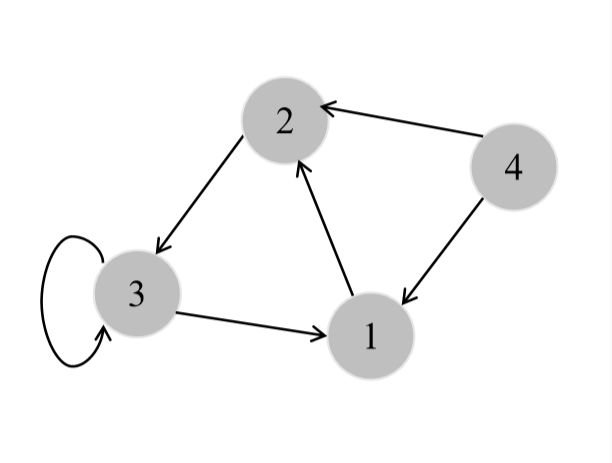}} 
\subfloat[$G_{\bar{A}_{\beta}}$]{\label{fig:0b}\includegraphics[width=0.32\textwidth]{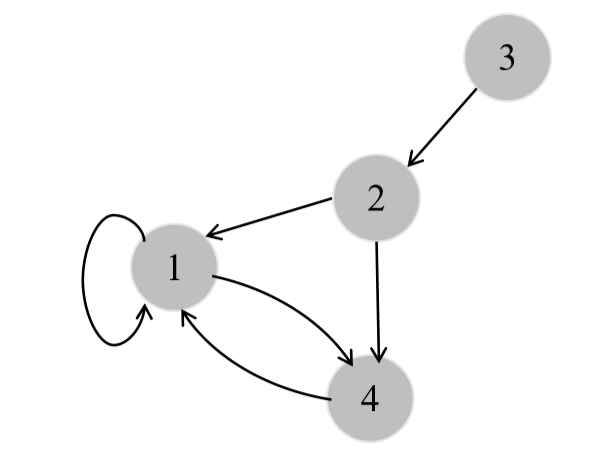}} 
\caption{Directed graph representations of the structured matrices $\bar A_{\alpha}$ and $\bar A_{\beta}$}
  \vspace{-1em}
\label{fig:0} 
\end{figure}

Consider the case $ n=4 $ with edge sets defined as $ E_{\alpha} = \{(1,2), (2,3), (3,1), (3,3)$, $(4,1),(4,2)\}$  and $ E_{\beta} = \{(1,1), (1,4),(2,1), (2,4),(3,2), (4,1)\} $. The zero patterns $ Z_{E_{\alpha}} $ and $ Z_{E_{\beta}} $ describe the system structure with matrices:
\begin{equation}
\bar{A}_{\alpha}= \begin{bmatrix}
    0 & \ast & 0 & 0\\
    0 & 0 & \ast & 0\\
    \ast & 0 & \ast & 0\\
    \ast & \ast & 0 & 0
\end{bmatrix},\quad
\bar{A}_{\beta}=\begin{bmatrix}
    \ast & 0 & 0 & \ast\\
    \ast & 0 & 0 & \ast\\
    0 & \ast & 0 & 0\\
    \ast & 0 & 0 & 0
\end{bmatrix}
\end{equation}
where $ \ast $ entries correspond precisely to the edges in $E$. The corresponding directed graphs with edge sets $E_{\alpha}$ and $E_{\beta}$, given by $G_{\bar A_{\alpha}}$ and $G_{\bar A_{\beta}}$ respectively, are shown in Figs. \ref{fig:0} (a) and (b).

For example, by assigning
\(
(\alpha_1,\ldots,\alpha_6)=(1,-2,3,4,-1,2)
\)
to the free entries of \(\bar A_{\alpha}\), we obtain a numerical
realization accordingly 
\[
A_{\alpha}=
\begin{bmatrix}
0 & 1 & 0 & 0\\
0 & 0 & -2 & 0\\
3 & 0 & 4 & 0\\
-1 & 2 & 0 & 0
\end{bmatrix}.
\]
Let $G = (V, E)$ be a directed graph (digraph) with vertex set $V = \{1, \dots, n\}$ and edge set $E \subseteq V \times V$, where $(i, j) \in E$ denotes a directed edge from vertex $i$ to vertex $j$.  An edge of the form $ (i, i) $ is called a \emph{loop} or \emph{self-loop}, indicating that its start and end vertices are identical.
A subset $ I \subseteq V $ is called an \emph{independent set} if there are no edges between any two vertices in $ I $, i.e.,
$
\forall u, v \in I \Rightarrow (u, v) \notin E.
$
 A \emph{path} is a sequence of edges connecting a sequence of vertices in a graph, where neither vertices nor edges are repeated. The length of a path is the number of edges traversed. A \emph{cycle} is a closed path, i.e., its start and end vertices are the same. \emph{Disjoint cycles} are cycles whose vertex sets do not intersect. A \emph{Hamiltonian cycle} is a cycle in a graph that visits every vertex exactly once and returns to the starting vertex. A $ k\text{-decomposition} $ of a graph $ G $ is a decomposition where the union of disjoint cycles covers exactly $ k $ vertices in $ G $. A decomposition covering all $ n $ vertices of $ G $ is called a \emph{Hamiltonian decomposition}. A graph $ G(V_1, V_2, E) $ is a \emph{bipartite graph} if and only if its vertex set $ V $ can be partitioned into two disjoint subsets $ V_1 $ and $ V_2 $ such that every edge $ e = (u, v) \in E $ satisfies $ u \in V_1 $ and $ v \in V_2 $, and no edges exist between vertices within the same subset in a bipartite graph. For a vertex subset $I \subseteq V_1$ ($I\subseteq V_2$), the \emph{neighborhood set} $ N(I) \subseteq V_2 $ ($N(I) \subseteq V_1$) satisfies:
$
\forall  v \in N(I)$, if there exists $u\in I$ such that $(u, v) \in E$ ($(v,u)\in E$).
In a bipartite graph $ G(V_1, V_2, E) $, a \emph{matching} is a subset of edges from $E$ where no two edges share a common vertex. A \emph{maximum matching} is a matching with the largest number of edges.
A \emph{perfect matching} of $G(V_1,V_2,E)$ is a matching covering all vertices of  $G(V_1,V_2,E)$. A bipartite graph has a perfect matching only if $|V_1| = |V_2|$.

{A digraph can be represented by a bipartite graph whose two vertex sets each contain a full copy of all vertices from the original graph. Specifically, for a digraph $G=(V,E)$, we define its corresponding bipartite graph as $\mathcal{B}(G)\doteq (V_{L},V_{R},E_{VV})$, where $V_{L}=V_{R}=V$ and $E_{VV}=E$ (for ease of notation, vertices in different partitions are allowed to share the same original label). Under this representation, any matching $\mathcal{M}$ in $\mathcal{B}(G)$ corresponds to a subgraph $G^{\cal M}\doteq (V,\mathcal{M})$ of $G$ with vertex set $V$ and edge set $\mathcal{M}\subseteq E_{VV}$.
With a slight abuse of terminology, ${\mathcal M}$ is also called a matching of $G$ and $G^{\cal M}$ is a decomposition of $G$ (by ${\cal M}$). From \cite[Lemma 3]{PequitoFramework2016}, $G^{\cal M}$ is a collection of disjoint cycles and paths.}

{In the bipartite graph $\mathcal{B}(G)=(V_{L}, V_{R}, E_{VV})$, if there exists a matching $\mathcal{M}$ of size $k$ that covers a left vertex set $S_L \subseteq V_L$ and a right vertex set $S_R \subseteq V_R$ with $|S_L| = |S_R| = k$, and if there exists a subset $S \subseteq V$ of size $k$ such that the copy of $S$ in the left part coincides with $S_L$ and its copy in the right part coincides with $S_R$, then we call $\mathcal{M}$ a (size-$k$) \emph{consistent matching}. If no such $S$ exists, we call $\mathcal{M}$ a (size-$k$) \emph{non-consistent matching}. For a consistent matching $\cal M$, $G^{\cal M}$ is a collection of disjoint cycles \cite{PequitoFramework2016}.   A perfect matching of ${\mathcal B}(G)$ is a size-$|V|$ consistent matching.}

\begin{figure}[htbp]
	\centering
	\subfloat[]{\label{fig:fa}\includegraphics[width=0.27\textwidth]{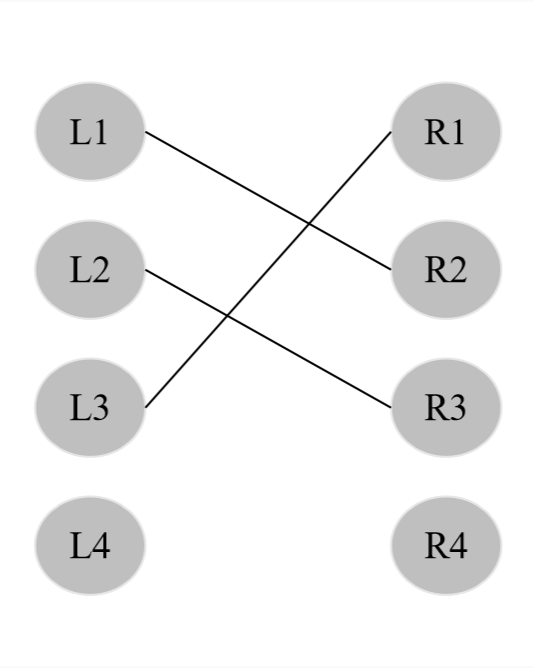}}
    \hspace{1cm}
	\subfloat[]{\label{fig:fb}\includegraphics[width=0.27\textwidth]{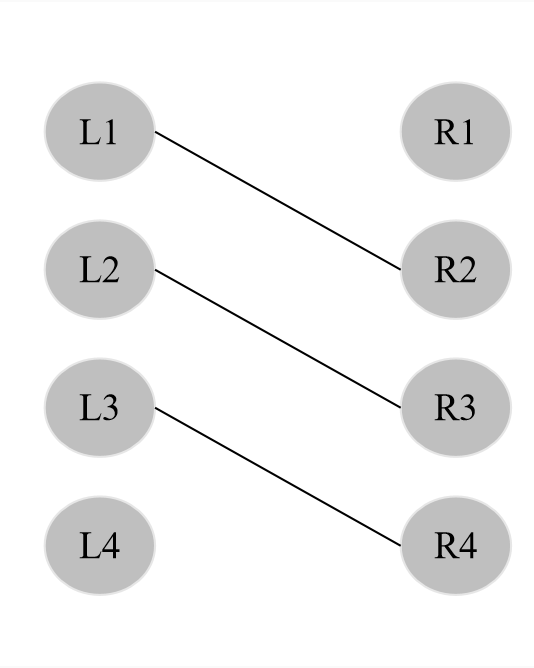}}
	\caption{Consistent and non-consistent matchings in a bipartite graph: (a) size-$3$ consistent matching, (b) size-$3$ non-consistent matching.}
\end{figure}

\begin{figure}[htbp]
	\centering
	\subfloat[]{\label{fig:dfa}\includegraphics[width=0.3\textwidth]{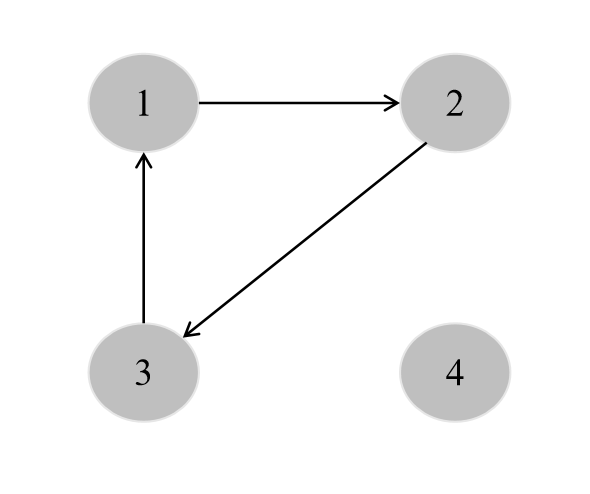}}
    \hspace{1cm}
	\subfloat[]{\label{fig:dfb}\includegraphics[width=0.3\textwidth]{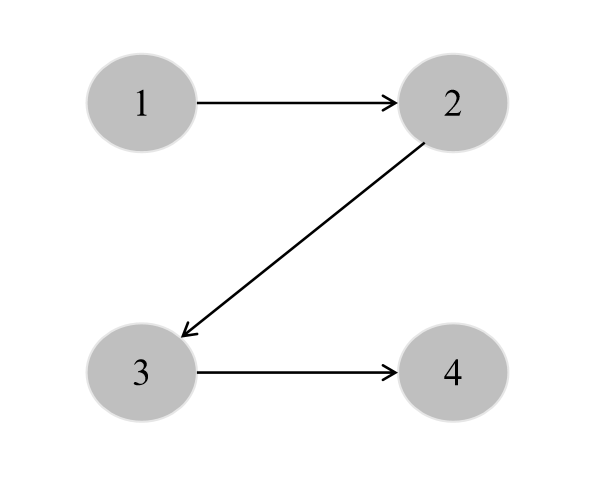}}
	\caption{Consistent and non-consistent matchings in the digraph: (a) size-$3$ consistent matching, (b) size-$3$ non-consistent matching. }
\end{figure}

{Fig.~\ref{fig:fa} shows a size-$3$ consistent matching, while Fig.~\ref{fig:fb} shows a size-$3$ non-consistent matching. Furthermore, by mapping Fig.~\ref{fig:fa} and Fig.~\ref{fig:fb} back to their corresponding digraphs,  a size-$3$ consistent matching  corresponds to a directed cycle of length $3$ (see Fig.~\ref{fig:dfa}), whereas a size-$3$ non-consistent matching corresponds to a directed path of length $3$ (see Fig.~\ref{fig:dfb}).}

\begin{definition}
The \emph{structural rank} (\emph{generic rank}) of a matrix $\bar{A} \in \{0,*\}^{n\times n}$, given by ${\rm grank}\, \bar A$ is defined as the maximum rank among all realizations of $\bar{A}$. Formally, it is the rank of the sparsity pattern of $\bar{A}$, assuming that all nonzero entries are algebraically independent.
\end{definition}

It has been shown in \cite{zhang2025generic} that given $\bar A\in \{0,*\}^{n\times n}$,
diagonalizability is a generic property in the sense that either almost all realizations of $\bar A$ are diagonalizable, or almost all of them are non-diagonalizable. This leads to the definition of structural diagonalizability.

\begin{definition} \label{def-diagonal}
    A structured matrix $\bar{A} \in \{0,*\}^{n \times n}$ is defined to be \emph{structurally diagonalizable} if, for almost all real value assignments to its free parameters, the resulting matrix is diagonalizable.
\end{definition}

\begin{lemma}[taken from \cite{zhang2025generic}]\label{lem:structural diagonalizability}
A structured matrix $\bar{A} \in \{0, \ast\}^{n \times n}$ is structurally diagonalizable if and only if the following equivalent conditions hold:
\begin{enumerate}
    \item $\mathrm{grank}(\bar{A}) = v(\bar{A})$, where $v(\bar{A})$ denotes the maximum number of vertices covered by disjoint cycles in the associated digraph ${G_{\bar A}}$.
    \item There exists a maximum matching in ${\cal B}({G_{\bar A}})$ such that the graph ${G_{\bar A}}$ decomposes into disjoint cycles and isolated vertices.
\end{enumerate}
\end{lemma}

According to \cref{lem:structural diagonalizability}, the structural diagonalizability of a structured matrix is determined by specific graph-theoretic properties of its associated digraph. Hence, the terms structural diagonalizability of a digraph and structural diagonalizability of its corresponding structured matrix can be used equivalently in this context. 

\begin{figure}[htbp]
	\centering
	\subfloat[\({\cal B}(G_{\bar{A}_{{\alpha}}})\)]{\includegraphics[width=0.27\textwidth]{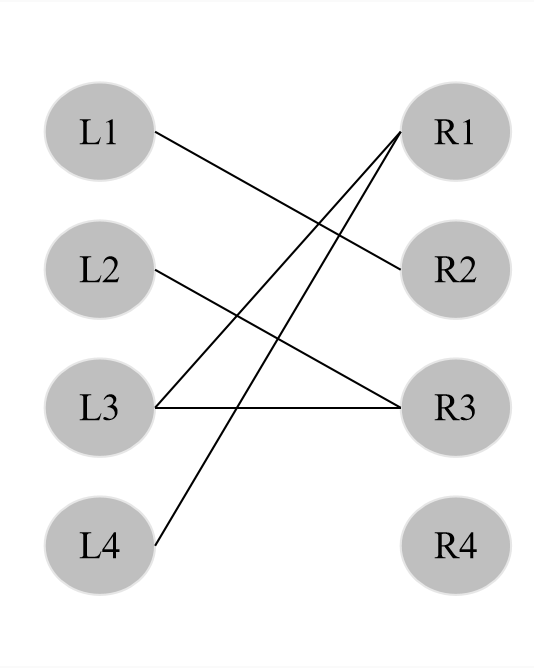}}
    \hspace{1cm}
	\subfloat[\({\cal B}(G_{\bar{A}_{{\beta}}})\)]{\includegraphics[width=0.27\textwidth]{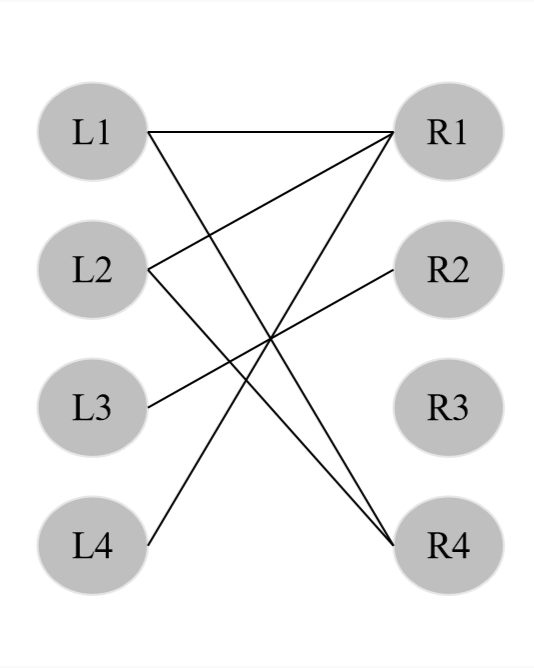}}
	\caption{The bipartite representations of \(G_{\bar{A}_{{\alpha}}}\) and \(G_{\bar{A}_{{\beta}}}\)}
    \label{fig:3}
\end{figure}

{The bipartite representations of \(G_{\bar{A}_{{\alpha}}}\) and \(G_{\bar{A}_{{\beta}}}\) from Fig.~\ref{fig:0} are illustrated in Fig.~\ref{fig:3}.  As shown in Fig.~\ref{fig:3}, the maximum matching size for both \({\cal B}(G_{\bar{A}_{{\alpha}}})\) and \({\cal B}(G_{\bar{A}_{{\beta}}})\) is 3. In \({\cal B}(G_{\bar{A}_{{\alpha}}})\), there exists a maximum matching that covers the vertex set \(\{1, 2, 3\}\).  And it is easy to see that the maximum matching that covers the vertex set \(\{1, 2, 3\}\) is a size-$3$ consistent matching, corresponding to the formation of cycle $1\rightarrow2\rightarrow3\rightarrow1$ in the original digraph. Additionally, the remaining vertex $4$ is isolated. In contrast, all maximum matchings in \({\cal B}(G_{\bar{A}_{{\beta}}})\) are size-$3$ non-consistent matchings, which correspond to chain-like structures in the digraph. From \cref{lem:structural diagonalizability}, \(G_{\bar{A}_{{\alpha}}}\) is structurally diagonalizable, whereas \(G_{\bar{A}_{{\beta}}}\) is not.}

All results in the following sections are stated in asymptotic form as \( n \to \infty \), using standard asymptotic notation. In particular, for functions \( f(n) \) and \( g(n) \):
\begin{itemize}
    \item \( f(n) = O(g(n)) \) means that there exist constants \( c > 0 \) and \( n_0 \) such that \( |f(n)| \leq c \, |g(n)| \) for all \( n > n_0 \).
    \item \( f(n) \in o(g(n)) \) means for any constant \( c > 0 \), there exists \( n_0 \) such that \( |f(n)| < c \, |g(n)| \) for all \( n > n_0 \).
    \item \( f(n) \in \omega(g(n)) \) means that for any constant \( c > 0 \), there exists \( n_0 \) such that \( |f(n)| > c \, |g(n)| \) for all \( n > n_0 \).
    \item $f(n) \sim g(n)$ means that $\lim_{n \to \infty} \frac{f(n)}{g(n)} = 1$.
\end{itemize}

In particular, \( f(n) \in o(1) \) means that \( \lim_{n \to \infty} f(n) = 0 \), and \( f(n) \in \omega(1) \) means that \( \lim_{n \to \infty} f(n) = +\infty \). The notation \( f(x) \ll g(x) \) (Vinogradov symbol) means that for all sufficiently large \( x \), \( f(x) < c \, g(x) \) holds for some constant \( c > 0 \).  ${\rm P}(\omega)$ denotes the probability of a specific outcome $\omega$ occurring in the probability space, while ${\mathbb E}X$ represents the expected value of a random variable $X$. Hereafter, \(\exp(\cdot)\) denotes the exponential function, i.e., \(\exp(x)=\text{e}^x\).

Using the random graph model enables us to study the probability distribution of the diagonalizability under different connection probabilities $p$. In the asymptotic regime where the edge probability is given by \( p(n) = \frac{\log n + c + o(1)}{n} \), the random graph is classified as {\emph{dense}} if \( c \in \omega(1) \), medium if \( c \) is a {\emph{constant}} (with respect to \( n \)), and {\emph{sparse}} if \( c \in -\omega(1) \).

The set of structurally diagonalizable digraphs with $n$ vertices is denoted by $\mathcal{D}^{n}$. It is formally defined as:
\[
\mathcal{D}^{n} = \left\{ G \mid G=(V,E) \text{ is structurally diagonalizable},\lvert V\rvert=n  \right\}.
\]

{To simplify the notation, we denote the event $G \in \mathcal{G}(n, p(n))$ by $G_{n,p(n)}$, and refer to $G_{n,p(n)}$ as a random digraph sampled from the random digraph model $\mathcal{G}(n, p(n))$. Similarly, we denote the event $G \in \mathcal{G}(n, p(n),q(n))$ by $G_{n,p(n)}^{q(n)}$, and refer to $G_{n,p(n)}^{q(n)}$ as a random digraph sampled from the random digraph model $\mathcal{G}(n, p(n),q(n))$.} For ease of reference, Table~\ref{tab:graph_classes} summarizes the main graph classes used throughout the paper. 
The table is intended as a notation guide; the formal definitions and characterizations are given in the corresponding sections indicated in the last column.

\begin{table}[t]
\centering
\caption{Summary of the main graph classes and random graph notations used throughout the paper.}
\label{tab:graph_classes}
\renewcommand{\arraystretch}{1.18}
\footnotesize

\begin{tabular}{c p{0.58\linewidth} p{0.2\linewidth}}
\hline
Notation & \multicolumn{1}{c}{Definition} & \multicolumn{1}{c}{Location} \\
\hline
\(\mathcal{G}(n,p,q)\)
&
The random digraph model on \(n\) labeled vertices in which each non-self-loop directed edge appears independently with probability \(p\), and each self-loop appears independently with probability \(q\).
&
Sec.~1 \\

\(\mathcal{G}(n,p)\)
&
The special case of \(\mathcal{G}(n,p,q)\) with \(q=p\).
&
Sec.~1 \\

\(G_{n,p(n)}\)
&
A random digraph sampled from the model \(\mathcal{G}(n,p(n))\).
&
Sec.~2 \\

\(G_{n,p(n)}^{q(n)}\)
&
A random digraph sampled from \(\mathcal{G}(n,p(n),q(n))\).
&
Sec.~2 \\

\(\mathcal{D}^{n}\)
&
The set of all \(n\)-vertex structurally diagonalizable digraphs.
&
Sec.~2 \\

\(\Gamma^{n}\)
&
The set of all \(n\)-vertex digraphs that do not admit a Hamiltonian decomposition.
&
\cref{lem:neighbor-connection} \\

\(\mathcal{F}_k^n\)
&
The set of \(n\)-vertex digraphs whose bipartite representations contain a minimal Hall-type violating set \(I\) with \(|I|=k\) and \(|N(I)|=k-1\).
&
\cref{def:Fk} \\

\(\mathcal{L}_k^n\)
&
The set of \(n\)-vertex digraphs whose bipartite representations admit a maximum matching of size \(n-k\) that decomposes the digraph into disjoint cycles and isolated vertices.
&
\cref{thm:snD} \\

\(\mathcal{H}_k^n\)
&
The set of \(n\)-vertex digraphs whose bipartite representations contain a size-\((n-k)\) consistent matching, and the remaining \(k\) vertices on either side of the bipartition are isolated.
&
\cref{def:Hkn} \\
\hline
\end{tabular}
\end{table}

A graph property is said to be monotone if adding edges to the original graph does not destroy the property \cite{frieze2015introduction}. Clearly, structural diagonalizability is not a monotone graph property, as adding a non-self-loop edge to a structurally diagonalizable empty graph makes it not structurally diagonalizable. This makes analyzing its asymptotic distribution in random graphs more challenging. 

\section{Graph-theoretic criteria for structural diagonalizability}
We first identify some easily verifiable sufficient or necessary graph-theoretic conditions for a digraph to be structurally diagonalizable. These characterizations will be helpful in deriving lower or upper bounds for the asymptotic probability of the diagonalizability of digraph structures.

\begin{lemma}[Sufficient condition for structural diagonalizability]\label{thm:Sufficient Conditions for structural diagonalizability}
A digraph is structurally diagonalizable if it has a Hamiltonian decomposition.
\end{lemma}
\begin{proof}
If a digraph admits a Hamiltonian decomposition, then it consists of disjoint cycles. This configuration satisfies the equivalent condition in \cref{lem:structural diagonalizability}, and thus the digraph is structurally diagonalizable.
\end{proof}

\begin{lemma}[Necessary and sufficient condition for the existence of a Hamiltonian decomposition, \cite{reinschke1988multivariable}]\label{lem:hd}
A digraph $G$ has a Hamiltonian decomposition if and only if its associated bipartite graph $\mathcal{B}(G)$ has a perfect matching.
  \end{lemma}

\begin{lemma}[Equivalent conditions for the nonexistence of a Hamiltonian decomposition]\label{lem:Equivalent Conditions for the Nonexistence of a Hamiltonian Decomposition}In a digraph $G$, a Hamiltonian decomposition does not exist if and only if there exists a subset \( I \subseteq V \) such that \( |N(I)| = |I| - 1 \), where \( V \) denotes one of the partite sets (left or right) in the corresponding bipartite $\mathcal{B}(G)$.
\end{lemma}
\begin{proof}
Using Hall's marriage theorem \cite{Hall1935} and \cref{lem:hd}: A Hamiltonian decomposition is absent if and only if $\mathcal{B}(G)$ has a subset $  I \subseteq V$(where \( V \) denotes one of the partite sets in $\mathcal{B}(G)$) such that $|N(I)| < |I|$. Let $I'$ be the minimal vertex subset satisfying $|N(I)| < |I|$. If $|N(I')| < |I'| - 1$, then removing any $v \in I'$ gives $|N(I' \setminus v)| \leq |N(I')| < |I'| - 1 = |I' \setminus v|$, contradicting the minimality of $I'$. Therefore, the minimal violating subset must satisfy $|N(I)| = |I| - 1$.

Conversely, if $\exists I \subseteq V$ with $|N(I)| = |I| - 1$, Hall's theorem implies no perfect matching. According to \cref{lem:hd}, we conclude that $G$ has no Hamiltonian decomposition.
\end{proof}

\begin{definition}\label{def:Fk}
Define $\mathcal{F}_k^n$ as the set of $n$-vertex digraphs whose bipartite graphs $\mathcal{B}=(V,V,E_{VV})$ satisfy:
\begin{enumerate}
    \item $\exists I \subseteq V$ with $|N(I)| = k - 1$ and $|I| = k$, where $V$ is the left or right vertex set of~$\mathcal{B}$;
    \item $\forall J \subseteq V$, if $|J| < k$, then $|N(J)| \neq |J| - 1$, where $V$ is the left or right vertex set of $\mathcal{B}$.
\end{enumerate}
\end{definition}
\begin{lemma}\label{lem:neighbor-connection}
$\mathcal{F}_k^n$ are pairwise disjoint and {$\bigsqcup_{k=1}^{\lceil n/2\rceil} \mathcal{F}_k^n = \Gamma^n$}, where $\Gamma^n$ is the set of all $n$-vertex digraphs without Hamiltonian decompositions, and $\lceil n/2 \rceil$ denotes the ceiling of $n/2$, i.e., the smallest integer greater than or equal to $n/2$.
\end{lemma}
\begin{proof}
By definition, $\mathcal{F}_k^n$ are disjoint. By \cref{def:Fk} and \cref{lem:Equivalent Conditions for the Nonexistence of a Hamiltonian Decomposition}, $\mathcal{F}_k^n \subseteq \Gamma^n$. Hall's marriage theorem states that a bipartite graph ${\cal B}=(V_L,V_R,E)$ has a perfect matching if and only if $ \forall I \subseteq V$, ($V=V_L$ or $V_R$), $|N(I)| \geq |I|$. It can be shown that ${\cal B}$ has a perfect matching if and only if
\begin{equation} \label{cond1}
    \begin{array}{c}
       \forall S \subseteq V_{L}, |S| \leq \frac{n}{2}, |N(S)| \geq |S|; \\
         \forall T \subseteq V_{R}, |T| \leq \frac{n}{2}, |N(T)| \geq |T|.
    \end{array}
\end{equation}

The only if direction is obvious. Now suppose condition \eqref{cond1} holds.
If $S \subseteq V_L$ with $|S| > \frac{n}{2}$ and $|N(S)| < |S|$, then setting $T = V_R \setminus N(S)$ yields $|N(T)|\le |V_L\backslash S|=|V_L|-|S|$ since $N(T)\subseteq V_L\backslash S$. As $|T|=|V_R|-|N(S)|> |V_R|-|S|$ and $|V_L|=|V_R|=n$, it follows that $|T|> |N(T)|$ and $|T|\le n-\frac{n}{2}$ (condition \eqref{cond1} yields $|N(S)|\ge \frac{n}{2}$), contradicting $\forall T \subseteq V_R$, $|T| \leq \frac{n}{2}$, $|N(T)| \geq |T|$. Hence, $\sum_{k=1}^{\lceil n/2\rceil} \mathcal{F}_k^n = \Gamma^n$.
\end{proof}

\begin{lemma}\label{lem:The neighborhood set satisfies the condition.}
In each digraph in $\mathcal{F}_k^n$, $k\ge 2$, every vertex in $N(I)$ (satisfying condition 1 of $\mathcal{F}_k^n$) connects to at least two vertices in $I$.
\end{lemma}
\begin{proof}
Select $I$ satisfying condition 1 of \cref{def:Fk}. If a vertex $v$ in $N(I)$ connects to only one vertex in $I$, removing $v$ yields $|N(I \setminus v)| = |I| - 1 - 1 = |I \setminus v| - 1$, contradicting condition 2 of $\mathcal{F}_k^n$.
\end{proof}

\begin{lemma}[Probability calculation for the non-existence of Hamiltonian decompositions in graphs]\label{thm:Probability Calculation for the Non-Existence of Hamiltonian Decompositions in Graphs}
For a n-vertex digraph $G$, we have
$\sum_{k=1}^{\lceil n/2\rceil} \mathrm{P}(G \in \mathcal{F}_k^n) = \mathrm{P}(G \in \Gamma^n)$.
\end{lemma}
\begin{proof}
From \cref{lem:neighbor-connection}, {$\bigsqcup_{k=1}^{\lceil n/2\rceil} \mathcal{F}_k^n = \Gamma^n$}, hence $\sum_{k=1}^{\lceil n/2\rceil} \mathrm{P}(G \in \mathcal{F}_k^n) = \mathrm{P}(G \in \Gamma^n)$.
\end{proof}

\begin{lemma}\label{thm:main-result}
Let $p(n) = \frac{\log n + c + o(1)}{n}$ and $0\leq q(n)\leq 1$. For $c \in \mathbb{R}$ or $c \in \pm \omega(1) $, we have
\begin{equation}
    \lim_{n\to\infty}\mathrm{P}(G_{n,p(n)}^{q(n)}\in \Gamma^n) = \lim_{n\to\infty}\mathrm{P}(G_{n,p(n)}^{q(n)} \in \mathcal{F}_1^n) + o(1).
\end{equation}
\end{lemma}

The proof can be found in Appendix A.  According to \cref{thm:main-result}, the probability that a digraph does not have a Hamiltonian decomposition asymptotically approaches the probability that the set $\mathcal{F}_1^n$ exists, i.e., the asymptotic probability that when viewing the digraph as a bipartite graph, there are isolated vertices on the left or right side. The approximation of this asymptotic probability helps us to further calculate the probability that a digraph has a Hamiltonian decomposition.

\begin{lemma}\label{lem:norn-1}
    If the bipartite graph ${\cal B}(G)$ corresponding to a digraph $G=(V,E)$ with $n$ vertices has a size-$n$ or size-$(n-1)$ consistent matching, then $G$ is structurally diagonalizable.
\end{lemma}

\begin{proof}
    If ${\cal B}(G)$ contains a size-$n$ consistent matching, then it possesses a perfect matching. By \cref{lem:hd}, this implies that $G$ admits a Hamiltonian decomposition, and hence it is structurally diagonalizable.

Now consider the case where ${\cal B}(G)$ contains a size-($n-1$) consistent matching ${\cal M}$. If a maximum matching of ${\cal B}(G)$ has size $n$,  it reduces to the previous case. If a maximum matching of ${\cal B}(G)$ has size $n-1$, then ${\cal M}$ is a maximum matching. Accordingly, the digraph $(V,{\cal M})$ consists of a collection of disjoint cycles and an isolated vertex. By Lemma \ref{lem:structural diagonalizability}, $G$ is structurally diagonalizable.
\end{proof}

From \cref{lem:norn-1}, when examining digraphs that are not structurally diagonalizable, it suffices to consider those for which the size of a consistent matching is at most $n-2$. 

{\
\begin{lemma}\label{thm:shangjie}
Consider a digraph $G$ with $n$ vertices. Let $u$ and $v$ be two distinct
original vertices ($u\ne v$). If all of the following three conditions hold
simultaneously:
\begin{enumerate}
    \item the left copy $u_L$ has no edge to any vertex in the right part of
    $\mathcal{B}(G)$;
    \item the right copy $v_R$ has no incident edge from any vertex in the
    left part of $\mathcal{B}(G)$;
    \item the cross-deleted bipartite subgraph
    $
        \mathcal{B}(G)\bigl[V_L\setminus\{u_L\},
        V_R\setminus\{v_R\}\bigr]
    $
    has a perfect matching of size $n-1$,
\end{enumerate}
then $G$ is not structurally diagonalizable.
\end{lemma}

\begin{proof}
First, observe that $u_L$ and $v_R$ are both isolated on their respective sides. This means the original vertex $u$ has no outgoing edge and the original vertex $v$ has no incoming edges. Therefore, neither $u$ nor $v$ can belong to a directed cycle. Since
$u\ne v$, it is clear that no consistent matching larger than $n-2$ exists.
Consequently, the digraph cannot contain a union of disjoint cycles covering
more than $n-2$ vertices.

However, due to Condition~3, we can construct a matching of size $n-1$.
Specifically, the perfect matching in the cross-deleted bipartite subgraph
matches every left vertex except $u_L$ and every right vertex except $v_R$,
and is therefore a matching of size $n-1$ in $\mathcal{B}(G)$.

Therefore, since $u_L$ is disconnected from all vertices on the opposite
side, the maximum matching size in $\mathcal{B}(G)$ is $n-1$. Nevertheless,
no maximum matching of this size corresponds to a structure of disjoint
cycles, since every consistent matching has size at most $n-2$. By
\cref{lem:structural diagonalizability}, $G$ is not structurally
diagonalizable.
\end{proof}
}

The three conditions in \cref{thm:shangjie} are required to hold simultaneously. Conditions~1 and~2 ensure that no consistent matching larger than size \(n-2\) exists, while Condition~3 guarantees the existence of a matching of size \(n-1\). Hence, the maximum matching size is strictly larger than the maximum consistent matching size, which implies non-diagonalizability. If any one of the three conditions is removed, this strict inequality may no longer hold.

\begin{figure}[htbp]
	\centering
	\subfloat[]{\label{fig:shangjie1}\includegraphics[width=0.25\textwidth]{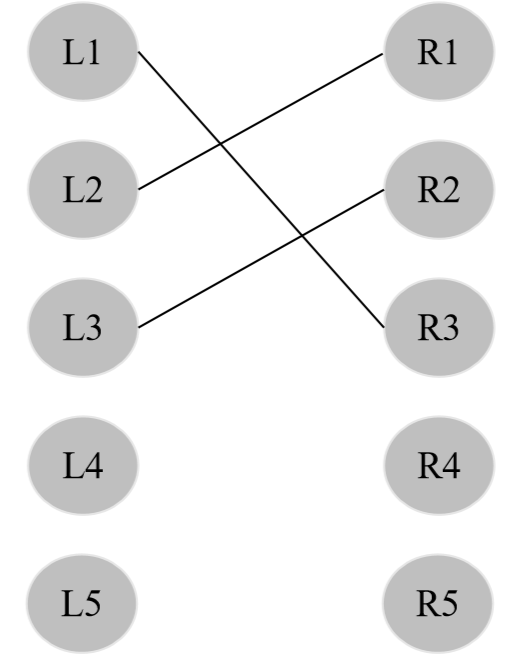}}
    \hspace{1cm}
	\subfloat[]{\label{fig:shangjie2}\includegraphics[width=0.25\textwidth]{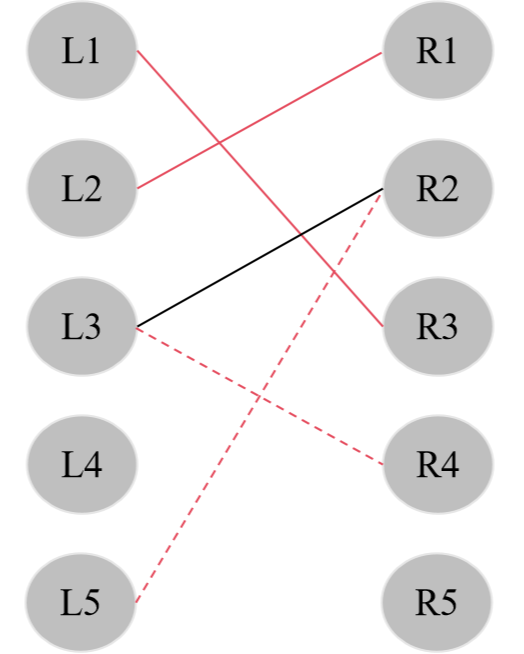}}
    \hspace{1cm}
    \subfloat[]{\label{fig:shangjie3}\includegraphics[width=0.25\textwidth]{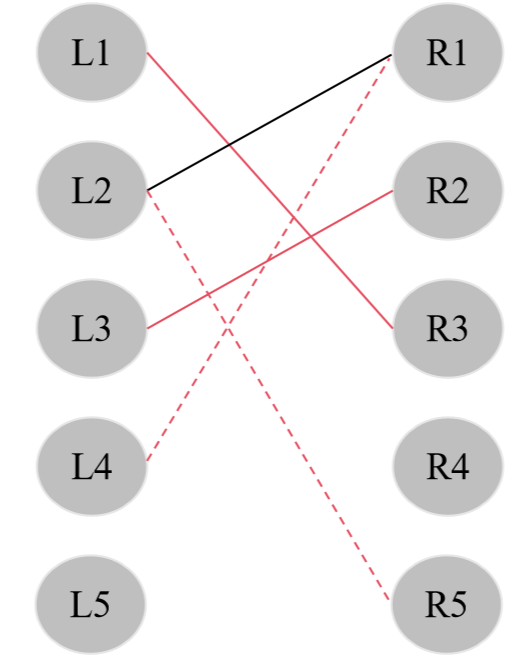}}
    \caption{Illustration of \cref{thm:shangjie}. (a) A size-\(3\) consistent matching \(\mathcal M\) in the bipartite representation of a \(5\)-vertex digraph. (b)--(c) Two cases in which all three conditions of \cref{thm:shangjie} hold simultaneously. Conditions~1 and~2 prevent the existence of a consistent matching larger than
size \(3\), {while the red solid and dashed edges together form the size-\(4\) perfect matching of the cross-deleted bipartite subgraph required by Condition~3.}} \label{fig:upper_bound_example}
  \end{figure}

Fig.~\ref{fig:shangjie1} shows a bipartite graph (for a 5-vertex digraph) with
a consistent matching of size 3. Based on this, Figs.~\ref{fig:shangjie2} and
\ref{fig:shangjie3} illustrate two cases that fulfill the conditions in
\cref{thm:shangjie}, {where the red solid and dashed lines together represent a
maximum matching of size 4. In particular, the red dashed lines indicate the
two matching edges used to match the remaining copies $u_R$ and $v_L$.
However, the enlarged matching is not a consistent matching. This violates
the necessary condition for structural diagonalizability, indicating that the
corresponding digraphs are not structurally diagonalizable.}

\begin{remark}
The configuration in \cref{thm:shangjie} is only one sufficient obstruction to structural diagonalizability. Additional configurations, such as consistent matchings of size no more than \(n-k\) for \(k\ge 3\) and larger non-consistent augmenting matchings, may yield tighter upper bounds. However, these
events generally overlap and require considerably more involved combinatorial analysis. We therefore focus on the tractable yet dominant configuration in \cref{thm:shangjie}, whose resulting bound is already relatively tight as shown in the numerical experiments. 
\end{remark}

\begin{theorem}\label{thm:snD} Let $\mathcal{L}_k^n$ denote the set of n-vertex digraphs $G$ whose bipartite graph admits a maximum matching of size $n-k$ that decomposes $G$ into disjoint cycles and isolated vertices. We have
\begin{equation*}
	\mathrm{P}\big(G_{n,p(n)} \in \mathcal{D}^{n}\big)
	=  \sum_{k=0}^n \mathrm{P}\big(G_{n,p(n)}\in \mathcal{L}_k^n\big).
\end{equation*}
Similarly,
 \begin{equation*}
	\mathrm{P}\big(G_{n,p(n)}^{q(n)} \in \mathcal{D}^{n}\big)
	=  \sum_{k=0}^n \mathrm{P}\big(G_{n,p(n)}^{q(n)} \in \mathcal{L}_k^n\big),
\end{equation*}
\end{theorem}
\begin{proof}
    It can be readily obtained from the graph-theoretic characterization of structural diagonalizability in \cref{lem:structural diagonalizability} and the disjointness of $\mathcal{L}_k^n$ (noting that every digraph has a unique maximum matching size).
\end{proof}

{
\begin{lemma} \label{thm:xiajie}
    For a $n$-vertex digraph $G$ and any integer $0 \leq k \leq n$, if its corresponding bipartite graph $\mathcal{B}(G)$ admits a consistent matching of size $n-k$, and the remaining $k$ vertices on one side of $\mathcal{B}(G)$ are isolated (i.e., they have no incident edges), then $G$ is structurally diagonalizable.
\end{lemma}

\begin{proof}
 Suppose the condition holds for some $k$: the bipartite graph $\mathcal{B}(G)$ admits a consistent matching $\cal M$ of size $n-k$, and the remaining $k$ vertices on one side are isolated. This implies that the original digraph $G$ contains a maximum matching $\cal M$ of size $n-k$. In addition, $\cal M$ decomposes $G$ into disjoint cycles and isolated vertices.  By the criterion for structural diagonalizability, $G$ is structurally diagonalizable.
\end{proof}}

Motivated by the above sufficient condition, we next identify a tractable subclass of structurally diagonalizable digraphs that will be used in the subsequent probabilistic analysis. 

\begin{definition}\label{def:Hkn}
Let \(\mathcal{H}_k^n\) denote the set of \(n\)-vertex digraphs whose bipartite representation contains a size-\((n-k)\) consistent matching, and the remaining \(k\) vertices on either side of the bipartition are isolated.
\end{definition}

Using the definition of \(\mathcal{L}_k^n\) from \cref{thm:snD}, we obtain $\mathcal{H}_{k}^n\subseteq\mathcal{L}_{k}^n$. Thus we have,
\begin{equation}\label{equ:xiajiegnp}
         \mathrm{P}\bigl(G_{n,p(n)} \in \mathcal{D}^{n}\bigr)
        \geq \sum_{k=0}^{n} \mathrm{P}\bigl(G_{n,p(n)} \in \mathcal{H}_k^n\bigr).
\end{equation}
Similarly,
\begin{equation}\label{equ:xiajiegnpq}
         \mathrm{P}\bigl(G_{n,p(n)}^{q(n)} \in \mathcal{D}^{n}\bigr)\geq  \sum_{k=0}^{n} \mathrm{P}\bigl(G_{n,p(n)}^{q(n)} \in \mathcal{H}_k^n\bigr).
    \end{equation}
    
\section{Asymptotic analysis of structural diagonalizability probability}

We primarily focus on two random graph models, namely, $\mathcal{G}(n,p)$ and $\mathcal{G}(n,p,q)$. We begin our investigation with the model $\mathcal{G}(n,p)$. Subsequently, we extend it to the $\mathcal{G}(n,p,q)$ model.

\subsection{Structural diagonalizability analysis in $\mathcal{G}(n,p)$}
We investigate the asymptotic probability of structural diagonalizability in digraphs under the $\mathcal{G}(n,p(n))$ model. We first set the edge probability $p(n)=\frac{\log n+c+o(1)}{n}$, and analyze how the limiting behavior of the constant $c$, ranging from positive infinity, convergence to a finite constant, to divergence to negative infinity, influences the asymptotic likelihood of diagonalizability across varying sparsity regimes.
\subsubsection{Structural diagonalizability analysis of dense graphs}
When the edge connection probability of a random digraph exceeds a certain threshold, the connectivity matrix of the digraph is highly likely to be structurally full rank, and the size of the maximum matching is highly likely to equal $n$. According to the first condition of \cref{lem:structural diagonalizability}, the random digraph is then highly likely to be structurally diagonalizable. Therefore, when studying the structural diagonalizability of dense random digraphs, we can consider the asymptotic probability of the digraph having a Hamiltonian decomposition ($\text{grank}(\bar{A}) = v(\bar{A}) = n$), in order to derive the asymptotic probability of structural diagonalizability for random digraphs.

  \begin{corollary}\label{cor:Asymptotic Probability} For the model ${\cal G}(n,p(n))$, we have
\begin{equation*}
    \lim_{n\to\infty}\mathrm{P}(G_{n,p(n)}\notin \Gamma^{n}) =
\begin{cases}
1 , & p(n)=\frac{\log n+c+o(1)}{n}\,\,,c\in\omega(1)\\
e^{-2e^{-c}} + o(1), & p(n)=\frac{\log n+c+o(1)}{n}\,\,,c \,\, \text{constant} \\
o(1), & p(n)=\frac{\log n+c+o(1)}{n}\,\,, c\in -\omega(1).
\end{cases}
\end{equation*}
\end{corollary}

\begin{proof}

According to \cref{thm:main-result}, the probability that a digraph does not have a Hamiltonian decomposition is asymptotically equivalent to the probability that the set $\mathcal{F}_1^n$ exists.

 Let $X_0$ denote the number of isolated vertices in the random bipartite graph $\mathcal{B}(G_{n,p(n)})$.
 Every possible edge between the two parts appears independently with probability $p(n)$.
A vertex is isolated if it has no edges to the opposite part. For a fixed vertex, there are $n$ potential edges to the opposite part, and the probability that none of these edges exist is $(1-p(n))^n$.
Since there are $2n$ vertices in $\mathcal{B}(G_{n,p(n)})$, by linearity of expectation, the expected number of isolated vertices is $\mathbb{E}X_0 = 2n(1-p(n))^n$.
Then, 
\begin{equation*}
    \lim_{n \to \infty} \mathbb{E}X_0 = \lim_{n \to \infty} 2n \cdot e^{-(\log n + c + o(1))} = 2e^{-c}.
\end{equation*}

Let $B_k^{(n)} = \mathbb{E}\binom{X_0}{k}$ denote the $k$th binomial moment of $X_0$. To compute $B_k^{(n)}$, fix $k$ and consider choosing $k$ distinct vertices from $\mathcal{B}(G_{n,p(n)})$. Let $a$ be the number of vertices selected from the left part and $b=k-a$ from the right part. The number of ways to select such a set is $\binom{n}{a}\binom{n}{k-a}$.
For these $k$ vertices to be isolated, all edges incident to them must be absent. The total number of potential edges incident to these vertices is $nk$, but the $ab$ edges between the chosen vertices in different parts are counted twice. Thus, exactly $nk - ab$ distinct edges must be absent, each independently with probability $1-p(n)$. Hence,
\[
\mathrm{P}(\text{all } k \text{ chosen vertices are isolated}) = (1-p(n))^{nk - ab}.
\]
Summing over all possible choices of the $k$ vertices,
\[
B_k^{(n)} = \sum_{a=0}^{k} \binom{n}{a}\binom{n}{k-a} (1-p(n))^{nk - ab}.
\]

For fixed $k$ as $n \to \infty$, we use the approximations:
\begin{align*}
\binom{n}{a} &\sim \frac{n^a}{a!}, \quad \binom{n}{k-a} \sim \frac{n^{k-a}}{(k-a)!}, \\
(1-p)^{nk - ab} &= \exp\!\bigl((nk - ab)\log(1-p)\bigr) \\
                &\sim \exp\!\bigl(-(nk - ab)p\bigr) \\
                &\sim  \exp\!\bigl(-k(\log n + c)\bigr) = n^{-k} e^{-kc}.
\end{align*}
Substituting into the expression for $B_k^{(n)}$,
\[
B_k^{(n)} \sim \sum_{a=0}^{k} \frac{n^a}{a!} \frac{n^{k-a}}{(k-a)!} n^{-k} e^{-kc}
          = e^{-kc} \sum_{a=0}^{k} \frac{1}{a!(k-a)!}.
\]
Using the identity $\sum_{a=0}^{k} \frac{1}{a!(k-a)!} = \frac{2^k}{k!}$, we obtain
\(
B_k^{(n)} \sim \frac{(2e^{-c})^k}{k!}.
\)

The binomial moments of a Poisson distribution $\operatorname{Po}(\lambda)$ are $\lambda^k/k!$. Since we have shown
\[
\lim_{n \to \infty} B_k^{(n)} = \frac{(2e^{-c})^k}{k!}
\]
for every fixed $k$, and since the Poisson distribution is uniquely determined by its moments, the method of moments \cite[Theorem 22.11]{frieze2015introduction} implies $X_0$ approximates a Poisson distribution with expectation $2e^{-c}$. The above derivation actually reflects the fact that the sum of independent binomial indicators turns a Poisson random variable in the many trials and finite expectation setting \cite[Theorem 3.1]{frieze2015introduction}.

Thus, we obtain $\lim_{n\to\infty}\mathrm{P}(X_0=0)=e^{-2e^{-c}}$. According to \cref{def:Fk}, we conclude that $\lim_{n\to\infty }\mathrm{P}(\mathcal{F}_{1}^n)=1-\lim_{n\to\infty}\mathrm{P}(X_0=0)=1-e^{-2e^{-c}}$.

From \cref{thm:main-result}, we have
\begin{equation*}
    \lim_{n \to \infty}\mathrm{P}(G_{n,p(n)} \in \Gamma^n) = \lim_{n\to\infty}\mathrm{P}(\mathcal{F}_1^n)+o(1)=1-e^{-2e^{-c}} + o(1),
\end{equation*}
and $\mathrm{P}(G_{n,p(n)}\notin \Gamma^{n}) = e^{-2e^{-c}} + o(1). $

Consequently, $\lim_{n \to \infty} \mathrm{P}(G_{n,p(n)}^{q(n)} \notin \Gamma^n) = e^{-2 e^{-c}}$ for the case $p(n) = \frac{\log n + c + o(1)}{n}$. To prove the case for $c\in \pm\omega(1)$, we can use monotonicity and the fact that $e^{-2e^{-c}} + o(1)\to 0$ if $c\in -\omega(1)$ and $e^{-2e^{-c}} + o(1)\to 1$ if $c\in \omega(1)$.
\end{proof}\par

According to \cref{lem:structural diagonalizability}, the existence of a Hamiltonian decomposition of a digraph is only a part of the structural diagonalizability of the digraph. Therefore, it is rough to estimate the lower bound of structural diagonalizability through Hamiltonian decomposition. For a dense digraph $G_{n,p(n)}$, when $p(n)=\frac{\log n+\omega(1)+o(1)}{n}$, we can obtain that the probability of the structural diagonalizability of a random digraph approaches $1$ asymptotically. Furthermore, we can conclude that when  $p(n)$ is an arbitrary constant (i.e., for dense graphs), the asymptotic probability of a random digraph structure being diagonalizable approaches $1$. The case \(p=0.5\) admits an additional combinatorial interpretation.
In the model \(\mathcal G(n,0.5)\), each of the \(n^2\) possible directed edges (including self-loops) is independently present with probability \(0.5\). Hence, every labeled digraph \(G\) on \(n\) vertices
has the same probability
\(
\mathrm P(G_{n,0.5}=G)
=
(0.5)^{|E(G)|}(0.5)^{n^2-|E(G)|}
=
2^{-n^2}.
\)
Since there are \(2^{n^2}\) labeled digraphs on \(n\) vertices when self-loops are allowed, we have
\(
\mathrm P(G_{n,0.5}\in\mathcal D^n)
=
\frac{|\mathcal D^n|}{2^{n^2}}.
\)
Thus, for \(p=0.5\), the probability
\(\mathrm P(G_{n,0.5}\in\mathcal D^n)\) is exactly the proportion of structurally diagonalizable digraphs among all \(n\)-vertex digraphs. Because this probability approaches $1$ as \(n\to\infty\), the proportion of such structurally diagonalizable digraphs also approaches $1$. This provides an affirmative answer to the conjecture proposed in~\cite{zhang2025generic}.

\subsubsection{Structural diagonalizability analysis of medium graphs}
In this part, we continue to investigate the asymptotic probability of structural diagonalizability for medium random digraphs.

\begin{theorem}\label{thm:gnpdown}
When $p(n)=\frac{\log n+c+o(1)}{n}$, we have
\begin{equation}
     \lim_{n\to\infty}\mathrm{P}(G_{n,p(n)}\in \mathcal{D}^{n}) \geq e^{-2e^{-c}}\big(1 + 2e^{-c}+e^{-2c}\big).
\end{equation}
\end{theorem}

\begin{proof}
    According to Eq.~\cref{equ:xiajiegnp} stemming from \cref{thm:xiajie}, we have
    \begin{align*}
        \lim_{n \to \infty} \mathrm{P}\bigl(G_{n,p(n)} \in \mathcal{D}^{n}\bigr)
        &\geq \lim_{n \to \infty} \sum_{k=0}^{n} \mathrm{P}\bigl(G_{n,p(n)} \in \mathcal{H}_k^n\bigr)\\
        &\geq \lim_{n \to \infty} \sum_{k=0}^{2} \mathrm{P}\bigl(G_{n,p(n)} \in \mathcal{H}_k^n\bigr).
    \end{align*}
    
      Furthermore, by \cref{cor:Asymptotic Probability}, as the number $n$ of vertices in the bipartite graph tends to infinity, the probability that the bipartite graph contains a perfect matching tends to $e^{-2e^{-c}}$. Consequently, as $n-k$ tends to infinity, the probability of having a size-(\(n-k\)) consistent matching also tends to $e^{-2e^{-c}}$. 
    For \(k=1\), we select one original vertex. The bipartite subgraph of the remaining \(n-1\) vertices has a size-\((n-1)\) consistent matching with asymptotic probability \(e^{-2e^{-c}}\). For the selected vertex, at least one of its two copies (left or right) must be isolated. By inclusion–exclusion, this probability \(2(1-p(n))^n - (1-p(n))^{2n-1}\). Multiplying by the number of choices \(\binom{n}{1}\) gives the contribution. The case \(k=2\) is handled similarly.  Therefore, we have
    \begin{align*}
        &\lim_{n \to \infty} \mathrm{P}\bigl(G_{n,p(n)} \in \mathcal{D}^{n}\bigr)\\
        &\geq \lim_{n \to \infty} \sum_{k=0}^{2} \mathrm{P}\bigl(G_{n,p(n)} \in \mathcal{H}_k^n\bigr)\\
        &= e^{-2e^{-c}} + \lim_{n\to\infty}\binom{n}{1} \times e^{-2e^{-c}} \times\big[ 2(1-p(n))^{n}-(1-p(n))^{2n-1}\big] \\
        &\quad + \lim_{n\to\infty}\binom{n}{2} \times e^{-2e^{-c}} \times \big[2(1-p(n))^{2n}-\big((1-p(n))^{2n-1}\big)^{2}\big]\\
        &= e^{-2e^{-c}} + \lim_{n\to\infty} n \times e^{-2e^{-c}} \times \big(2e^{-(\log n + c + o(1))}-e^{-2(\log n + c + o(1))}\big) \\
        &\quad + \lim_{n\to\infty} \frac{n(n-1)}{2} \times e^{-2e^{-c}} \times \big(2e^{-2(\log n + c + o(1))}-e^{-4(\log n + c + o(1))}\big)\\
        &= e^{-2e^{-c}} + \lim_{n\to\infty} n \times e^{-2e^{-c}} \times (2\frac{e^{-c}}{n}-\frac{e^{-2c}}{n^2}) \\
        &\quad + \lim_{n\to\infty} \frac{n(n-1)}{2} \times e^{-2e^{-c}} \times (2\frac{e^{-2c}}{n^2}-\frac{e^{-4c}}{n^4})\\
        &= e^{-2e^{-c}}\bigl(1 + 2e^{-c} + e^{-2c}\bigr).
    \end{align*}
\end{proof}

\begin{theorem}\label{thm:gnpup}
Let \(p(n)=\frac{\log n+c+o(1)}{n}\). Then,
\[
\lim_{n\to\infty}\mathrm{P}\bigl(G_{n,p(n)}\in\mathcal{D}^n\bigr)
\leq 1-e^{-2e^{-c}}e^{-2c}.
\]
\end{theorem}

\begin{proof}
It is straightforward to see that
\[
\mathrm{P}\bigl(G_{n,p(n)}\in\mathcal{D}^n\bigr)
=1-\mathrm{P}\bigl(G_{n,p(n)}\notin\mathcal{D}^n\bigr).
\]

By \cref{thm:shangjie}, a digraph fails to be structurally diagonalizable if
the following conditions hold for some pair of distinct original vertices
$u$ and $v$:
\begin{enumerate}
    \item The left copy $u_L$ has no edge to any vertex in the right part of
    the bipartite graph.
    \item The right copy $v_R$ has no incident edge from any vertex in the
    left part of the bipartite graph.
    \item After deleting $u_L$ and $v_R$, the remaining bipartite subgraph has
    a perfect matching of size $n-1$.
\end{enumerate}

The cross-deleted bipartite subgraph has $n-1$ vertices in each part. It follows from \cref{cor:Asymptotic Probability} that the probability
that this subgraph admits a perfect matching approaches
\(e^{-2e^{-c}}\). 

For a fixed ordered pair $(u_L,v_R)$, the first two conditions require exactly
$2n-1$ distinct edges to be absent, because the edge $(u_L,v_R)$ is counted
in both the row of $u_L$ and the column of $v_R$. These edges are disjoint
from those in the cross-deleted subgraph. Hence, by independence, all three
conditions hold with probability
$(1-p(n))^{2n-1}(e^{-2e^{-c}}+o(1))$.

There are \(n(n-1)\) choices of the ordered pair \((u_L,v_R)\).
These events are pairwise disjoint, since the perfect matching in the
cross-deleted subgraph guarantees that \(u_L\) and \(v_R\) are the unique
isolated left and right vertices, respectively. 
Thus, we have
\begin{align*}
  &\lim_{n\to\infty} \mathrm{P}\bigl(G_{n,p(n)} \in \mathcal{D}^n\bigr)\\
  &= 1 - \lim_{n\to\infty} \mathrm{P}\bigl(G_{n,p(n)} \notin \mathcal{D}^n\bigr)\\
&\leq
1-\lim_{n\to\infty}
n(n-1)(1-p(n))^{2n-1}
\bigl(e^{-2e^{-c}}+o(1)\bigr)\\
&=1-\lim_{n\to\infty}n(n-1)e^{-(2n-1)p(n)}\bigl(e^{-2e^{-c}}+o(1)\bigr)\\
&=1-\lim_{n\to\infty}n(n-1)e^{-\frac{(2n-1)(\log n+c+o(1))}{n}}\bigl(e^{-2e^{-c}}+o(1)\bigr)\\
&=1-\lim_{n\to\infty}n(n-1)e^{-2(\log n+c+o(1))}\bigl(e^{-2e^{-c}}+o(1)\bigr)\\
&=1-e^{-2e^{-c}}e^{-2c}.
\end{align*}
\end{proof}

From \cref{thm:gnpdown} and \cref{thm:gnpup}, when $p(n)=\frac{\log n+c+o(1)}{n}$, we can obtain that
\begin{equation}\label{twobounds}
    e^{-2e^{-c}}\big(1 + 2e^{-c}+e^{-2c}\big) \leq \lim_{n\to\infty}\mathrm{P}(G_{n,p(n)}\in  \mathcal{D}^{n})\leq 1-e^{-2e^{-c}} e^{-2c}.
\end{equation}

Moreover, it is straightforward to prove by differentiation that the upper bound is always greater than or equal to the lower bound.
  We have conducted an analysis of the structural diagonalizability for dense random digraphs and medium random digraphs. For the medium random digraphs, as can be seen from the upper and lower bounds, when $c\in -\omega(1)$, meaning the random digraph becomes sparse, the lower bound for structural diagonalizability tends to $0$ and the upper bound tends to $1$. This conclusion is trivial. Therefore, further analysis is still required for sparse random digraphs.
  \subsubsection{Structural diagonalizability analysis of sparse graphs}
Consider the edge probability $ p(n) = \frac{c}{n} $, where $ c $ is a positive constant.

An isolated edge refers to two distinct nodes $u$ and $v$ such that there is a directed edge from $u$ to $v$, and $u$ has out-degree $1$ and in-degree $0$, while $v$ has out-degree $0$ and in-degree~$1$.

{\begin{lemma}\label{lem:gulibian}
    A digraph $G$ containing an isolated edge is not structurally diagonalizable.
\end{lemma}}
\begin{proof}
    If $G$ has an isolated edge, then, by definition, this edge does not belong to any cycle, which means it cannot be contained by any disjoint cycle structure formed by a maximum matching. However, the two vertices connected by the isolated edge are not mutually independent. By Lemma \ref{lem:structural diagonalizability}, $G$ is not structurally diagonalizable.
\end{proof}

 From \cref{lem:gulibian}, a digraph containing an isolated edge is never structurally diagonalizable. Therefore, when computing an upper bound for $\mathrm{P}(G_{n,p(n)}\in\mathcal{D}^{n})$, we can approximate it by the asymptotic probability that $G_{n,p(n)}$ contains no isolated edges. This is exactly the idea of our proof for the following result, whose complete proof is given in Appendix~B.

\begin{theorem} \label{thm:SDB=0}
    If $ p(n) = \dfrac{c}{n} $ where $ c$ is a positive constant, we have
\begin{equation*}
\lim_{n \to \infty} \mathrm{P}(G_{n,p(n)}\in\mathcal{D}^{n})=0.
\end{equation*}
\end{theorem}

    The above theorem reveals that the asymptotic probability of structural diagonalizability for extremely sparse random digraphs is $0$.
So far, we have characterized the asymptotic probability distribution of structural diagonalizability for the $\mathcal{G}(n,p)$ random graph model across different edge probabilities $p$.

Notably, it can be verified that both lower and upper bounds in the medium graph region are monotonically increasing functions with respect to $c$. Overall, it seems that denser graphs are more likely to be structurally diagonalizable (verified by our simulations in \cref{sec-simulation}). This implies that even though structural diagonalizability is not a monotone graph property, it may possess this property in the statistical sense.

In what follows, we extend these analyses to the random network model ${\cal G}(n,p,q)$.

\subsection{ Structural diagonalizability analysis in $\mathcal{G}(n,p,q)$}
We continue with the random digraph model $\mathcal{G}(n, p(n), q(n))$. First, we extend \cref{cor:Asymptotic Probability} to this model. 
                                                  \begin{corollary}\label{cor:Asymptotic Probability gnpq}For the model $\mathcal{G}(n, p(n), q(n))$, $0\le q(n)\le 1$, we have 
\begin{equation*}
    \lim_{n\to\infty}\mathrm{P}(G_{n,p(n)}^{q(n)}\notin \Gamma^{n}) =
\begin{cases}
1 , & p(n)=\frac{\log n+c+o(1)}{n},\, c \in \omega(1)\\
e^{-2(1-q(n))e^{-c}}, & p(n)=\frac{\log n+c+o(1)}{n}\,\,, c \,\,\text{constant} \\
o(1), & p(n)=\frac{\log n+c+o(1)}{n}\,\,, c \in -\omega(1).
\end{cases}
\end{equation*}
\end{corollary}
\begin{proof}
  Let $p(n)=\frac{\log n+c+o(1)}{n}$ and $0\leq q(n)\leq 1$. The proof follows a similar idea to that of \cref{cor:Asymptotic Probability}, with the only difference being that the probability of a vertex forming a self-loop is changed from $p(n)$ to $q(n)$. Let $X_0$ denote the number of isolated vertices in the random bipartite graph $\mathcal{B}(G_{n,p(n)}^{q(n)})$.
 \begin{equation*}
\mathbb{E}X_0 = 2n\cdot(1-q(n))\cdot(1-p(n))^{n-1}.
\end{equation*}
Then, we have
\begin{equation*}
    \lim_{n \to \infty} \mathbb{E}X_0 = \lim_{n \to \infty} 2n \cdot(1-q(n))\cdot e^{-(\log n + c + o(1))} =2(1-q(n))\cdot e^{-c}.
\end{equation*}

Similar to the analysis of \cref{cor:Asymptotic Probability}, $X_0$ follows a Poisson distribution. We thus obtain $\lim_{n\to\infty}\mathrm{P}(X_0=0)=e^{-2(1-q(n)) e^{-c}}$. According to \cref{def:Fk}, we conclude that $\lim_{n\to\infty }\mathrm{P}(\mathcal{F}_{1}^n)=1-\lim_{n\to\infty}\mathrm{P}(X_0=0)=1-e^{-2(1-q(n))e^{-c}}$.

From \cref{thm:main-result}, we have
\begin{equation*}
    \lim_{n \to \infty}\mathrm{P}(G_{n,p(n)}^{q(n)} \in \Gamma^n) = \lim_{n\to\infty}\mathrm{P}(\mathcal{F}_1^n)+o(1)=1-e^{-2(1-q(n))e^{-c}} + o(1).
\end{equation*}
Thus, we can use monotonicity and the fact that $e^{-2(1-q(n))e^{-c}} + o(1)\to0$ if $c\in -\omega(1)$ and $e^{-2(1-q(n))e^{-c}} + o(1)\to 1$ if $c\in \omega(1)$.
\end{proof}\par

\begin{theorem}\label{thm:gpqdown}
When $p(n)=\frac{\log n+c+o(1)}{n}$ and $0\leq q(n)\leq 1$, we have
\begin{equation*}
     \lim_{n\to \infty}\mathrm{P}(G_{n,p(n)}^{q(n)}\in \mathcal{D}^n) \geq e^{-2(1-q(n))e^{-c}}[1 + 2(1-q(n))e^{-c}+(1-q(n))^{2}e^{-2c}].
\end{equation*}
\end{theorem}
\begin{proof}
This follows the same proof idea as \cref{thm:gnpdown}, with the probability of a vertex forming a self-loop changing from $p(n)$ to $q(n)$. The specific asymptotic probability calculations are given below:
\begin{align*}
    &\lim_{n \to \infty} \mathrm{P}\bigl(G_{n,p(n)}^{q(n)} \in \mathcal{D}^{n}\bigr)\\
    &\geq \lim_{n \to \infty} \sum_{k=0}^{2} \mathrm{P}\bigl(G_{n,p(n)}^{q(n)} \in \mathcal{H}_k^n\bigr)\\
    &= e^{-2(1-q(n))e^{-c}} + \lim_{n\to\infty}\binom{n}{1} \times e^{-2(1-q(n))e^{-c}} \times (1-q(n))\\
    &\quad\times \big[2(1-p(n))^{n-1}-(1-p(n))^{2n-2}\big] \quad + \lim_{n\to\infty}\binom{n}{2} \times e^{-2(1-q(n))e^{-c}} \\
    &\quad \times (1-q(n))^{2}\times \big[2(1-p(n))^{2n-2}-(1-p(n))^{4n-4}\big]\\
    &= e^{-2(1-q(n))e^{-c}}\big(1 + \lim_{n\to\infty} n  \times 2(1-q(n))e^{-(\log n + c + o(1))} \\
    &\quad + \lim_{n\to\infty} \frac{n(n-1)}{2} \times 2(1-q(n))^{2}e^{-2(\log n + c + o(1))} \big)\\
    &= e^{-2(1-q(n))e^{-c}}\bigl(1 + 2(1-q(n))e^{-c} + (1-q(n))^{2}e^{-2c}\big).
\end{align*}
\end{proof}

{
\begin{theorem}\label{thm:gpqup}
When $p(n)=\frac{\log n+c+o(1)}{n}$ and $0\leq q(n)\leq 1$, we have
\begin{equation*}
    \lim_{n\to \infty}\mathrm{P}(G_{n,p(n)}^{q(n)}\in  \mathcal{D}^{n})\leq 1-(1-q(n))^{2}\cdot e^{-2(1-q(n))e^{-c}} e^{-2c}.
\end{equation*}
\end{theorem}
\begin{proof}
We employ the same method of proof as used for \cref{thm:gnpup}. For a fixed
ordered pair $(u_L,v_R)$, the first two conditions in \cref{thm:shangjie} require
two self-loops and $2n-3$ distinct non-self-loop edges to be absent, with
probability $(1-q(n))^2(1-p(n))^{2n-3}$. The cross-deleted bipartite subgraph has $n-1$ vertices in each part. After
pairing $v_L$ with $u_R$, it differs from the bipartite representation of
$\mathcal G(n-1,p(n),q(n))$ only in the edge $(v_L,u_R)$, whose probability
is $p(n)$ rather than $q(n)$. Since the contribution of this exceptional edge
is $o(1)$, the same argument as in
\cref{cor:Asymptotic Probability gnpq} gives that the probability of a perfect
matching is
\(
e^{-2(1-q(n))e^{-c}}+o(1).
\)

There are
$n(n-1)$ ordered choices of $(u_L,v_R)$, and the corresponding
events are pairwise disjoint. The detailed calculations are as follows:
\begin{align*}
  &\lim_{n\to\infty}\mathrm P\bigl(G_{n,p(n)}^{q(n)}\in\mathcal D^n\bigr)\\
  &=1-\lim_{n\to\infty}\mathrm P\bigl(G_{n,p(n)}^{q(n)}\notin\mathcal D^n\bigr)\\
  &\leq 1-\lim_{n\to\infty}n(n-1)
  \bigl(e^{-2(1-q(n))e^{-c}}+o(1)\bigr)
  \cdot (1-q(n))^2(1-p(n))^{2n-3}\\
  &=1-\lim_{n\to\infty}n(n-1)(1-q(n))^2
  e^{\bigl(-(2n-3)p(n)+O(np(n)^2)\bigr)}
  \bigl(e^{-2(1-q(n))e^{-c}}+o(1)\bigr)\\
  &=1-\lim_{n\to\infty}n(n-1)(1-q(n))^2
  e^{-(2n-3)\frac{\log n+c+o(1)}{n}}
  \bigl(e^{-2(1-q(n))e^{-c}}+o(1)\bigr)\\
  &=1-\lim_{n\to\infty}n(n-1)(1-q(n))^2
  e^{-2(\log n+c+o(1))}
  \bigl(e^{-2(1-q(n))e^{-c}}+o(1)\bigr)\\
  &=1-(1-q(n))^2e^{-2(1-q(n))e^{-c}}e^{-2c}.
\end{align*}
\end{proof}
}


In the random graph model $\mathcal{G}(n,p(n),q(n))$, according to the upper and lower bounds on the diagonalizability probability in \cref{thm:gpqdown} and \cref{thm:gpqup}, we can conclude that when $p(n)=\frac{\log n+c+o(1)}{n}, 0\le q(n)\le 1$, 
\begin{align*}
    e^{-2(1-q(n))e^{-c}}\big(1 + &2(1-q(n))e^{-c}+(1-q(n))^{2}e^{-2c}\big) \leq \lim_{n\to \infty}\mathrm{P}(G_{n,p(n)}^{q(n)}\in  \mathcal{D}^{n})\\&\leq
    1-(1-q(n))^{2}\cdot e^{-2(1-q(n))e^{-c}} e^{-2c}.
\end{align*}It can be verified that when $q(n)=1$, both the lower and upper bounds above become $1$, which is consistent with the fact that graphs with self-loops on all nodes are structurally diagonalizable (c.f. \cref{lem:structural diagonalizability}). Moreover, when $q(n)\in o(1)$, the above bounds reduce to \cref{twobounds}, consistent with \cref{thm:gnpdown} and \cref{thm:gnpup}.

One can easily verify that the above bounds are consistent with both
special cases: setting
\(q(n)=p(n)=(\log n+c+o(1))/n\) recovers the corresponding bounds for
the \(\mathcal{G}(n,p)\) model, whereas setting \(q(n)=1\) makes both
bounds equal to \(1\), in agreement with the fact that the self-loops
form a cycle cover of all vertices.

To provide a clearer overview of the main probabilistic results, \cref{table:asymptotic_results} summarizes the asymptotic probability of structural diagonalizability for the random graph models \(\mathcal{G}(n,p)\) and \(\mathcal{G}(n,p,q)\) across different regimes.

\begin{table}[t]
\centering
\caption{Asymptotic probability of structural diagonalizability for the random graph models \(\mathcal{G}(n,p)\) and \(\mathcal{G}(n,p,q)\).}
\label{table:asymptotic_results}
\renewcommand{\arraystretch}{1.25}
\setlength{\tabcolsep}{3.2pt}
\footnotesize

\begin{tabular}{p{0.12\linewidth} p{0.14\linewidth} p{0.28\linewidth} p{0.35\linewidth}}
\toprule
\textbf{Model} 
& \textbf{Regime} 
& \textbf{Edge probability} 
& \textbf{Asymptotic probability} \\
\midrule

\multirow{3}{*}{\(\mathcal{G}(n,p)\)}
& Dense graph
& \(\begin{gathered}
p(n)=\dfrac{\log n+c+o(1)}{n}\\
c\in\omega(1)
\end{gathered}\)
& \(1\) \\

\cmidrule(lr){2-4}

& Medium graph
& \(\begin{gathered}
p(n)=\dfrac{\log n+c+o(1)}{n}\\
c\in\mathbb{R}
\end{gathered}\)
& \(\begin{gathered}
[\,e^{-2e^{-c}}(1+2e^{-c}+e^{-2c}),\\
1-e^{-2e^{-c}}e^{-2c}\,]
\end{gathered}\) \\

\cmidrule(lr){2-4}

& Sparse graph
& \(\begin{gathered}
p(n)=\dfrac{c}{n} \quad 
c>0
\end{gathered}\)
& \(0\) \\

\midrule

\multirow{2}{*}{\(\mathcal{G}(n,p,q)\)}
& Dense graph
& \(\begin{gathered}
p(n)=\dfrac{\log n+c+o(1)}{n}\\
c\in\omega(1),\ 0\le q(n)\le 1
\end{gathered}\)
& \(1\) \\

\cmidrule(lr){2-4}

& Medium graph
& \(\begin{gathered}
p(n)=\dfrac{\log n+c+o(1)}{n}\\
c\in\mathbb{R},\ 0\le q(n)\le 1
\end{gathered}\)
& \(\begin{gathered}
\bigl[
e^{-2(1-q(n))e^{-c}}
\bigl(1+2(1-q(n))
\\e^{-c}
 +(1-q(n))^{2}e^{-2c}\bigr),
1-\\(1-q(n))^{2}e^{-2c}
e^{-2(1-q(n))e^{-c}}
\bigr]
\end{gathered}\) \\

\bottomrule
\end{tabular}
\end{table}

\section{Experiments and validation}\label{sec-simulation}

In this section, we conduct numerical simulations to evaluate the probability of structural diagonalizability in random graphs. 
The MATLAB code  used to reproduce the numerical experiments is publicly available in the GitHub repository \url{https://github.com/ewbo123/diagonalizable-random-structures}.

\subsection{The simulation results of the $\mathcal{G}(n,p)$) random graph model}
The simulation framework contains two key parameters:
\begin{itemize}
    \item $p$: Connection probability of graphs, determined by $n$ and $c$
    \item $n$: Number of vertices in generated graphs
  \end{itemize}

 This experiment investigates the consistency between simulation results and theoretical bounds across different $c$ values under fixed network sizes. Two connection probability regimes are examined:

Case 1: $p = \frac{\log n + c+o(1)}{n}$
\begin{itemize}
    \item {Fixed $n =4000$ with $800$ samples per configuration}
    \item $c \in [-1, 2]$ sampled at 0.1 intervals for smooth curves
    \item Theoretical bounds: \begin{equation*}
    e^{-2e^{-c}}\big(1 + 2e^{-c}+e^{-2c}\big) \leq \lim_{n\to \infty}\mathrm{P}(G_{n,p}\in  \mathcal{D}^{n})\leq 1-e^{-2e^{-c}} e^{-2c}.
\end{equation*}
\end{itemize}

Case 2: $p = \frac{c}{n}$
\begin{itemize}
    \item Fixed $n = 4000$ with $800$ samples
    \item $c \in [1, 5]$ to test boundary effects
    \item Theoretical prediction: $\lim_{n\to \infty}\mathrm{P}(G_{n,p}\in  \mathcal{D}^{n}) \equiv 0$.
\end{itemize}

\begin{figure}[htbp]
\centering
  {\includegraphics[width=0.6\textwidth]{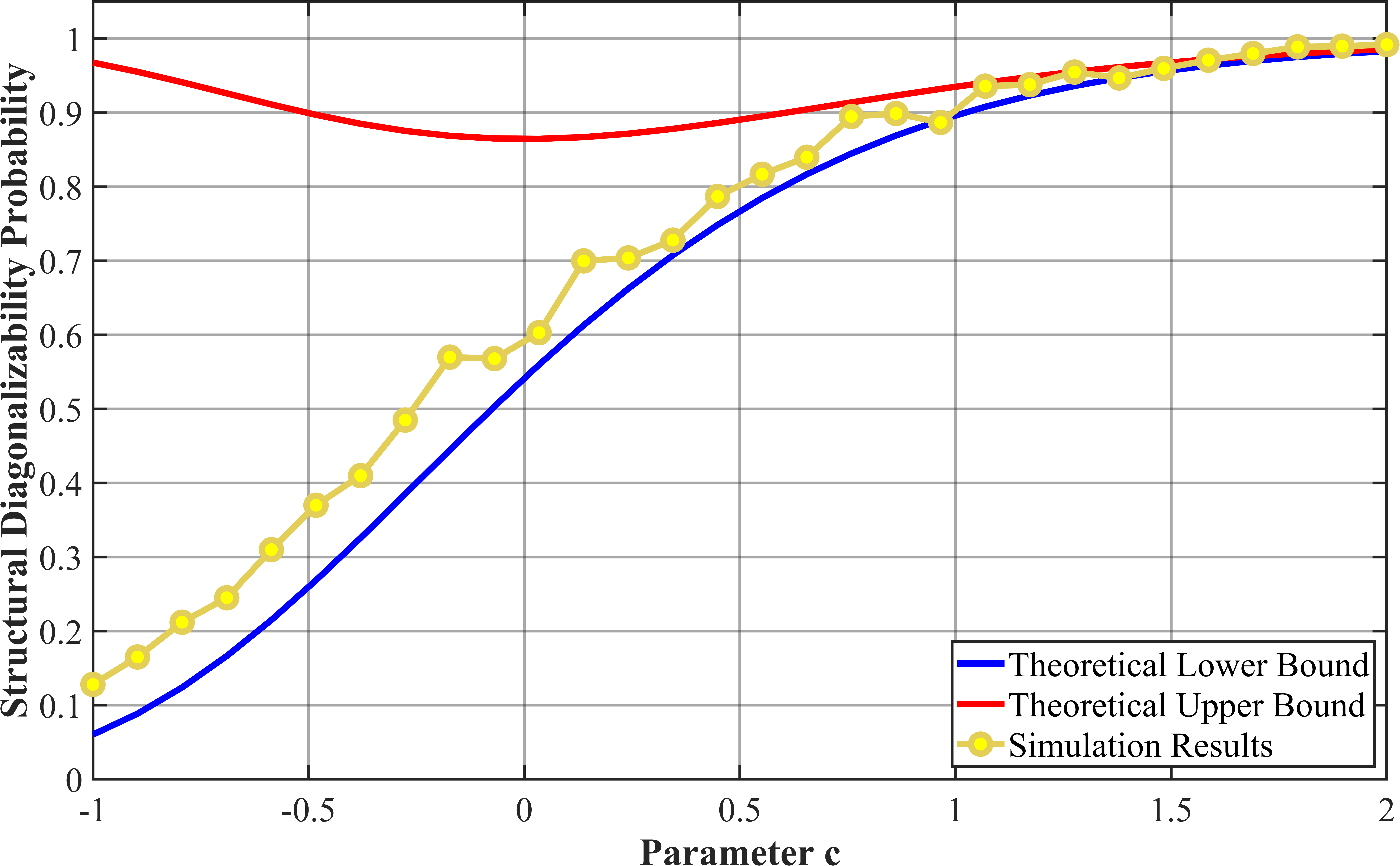}}
\caption{Comparison between simulation results and theoretical upper/lower bounds for $p=\frac{\log n+c+o(1)}{n}$ with fixed $n=4000$ and varying $c$}
\vspace{-1em}
\label{fig:log_n_n} 
\end{figure}

From Fig.~\ref{fig:log_n_n}, it can be observed that for random digraphs with vertex counts $n=4000$, where $800$ graph instances are generated per parameter configuration, the variation trend of the structural diagonalizability probability over $c \in [-1, 2]$ closely aligns with the theoretically derived upper and lower bounds.

\begin{figure}[htbp]  
\centering
{\label{fig:a}\includegraphics[width=0.6\textwidth]{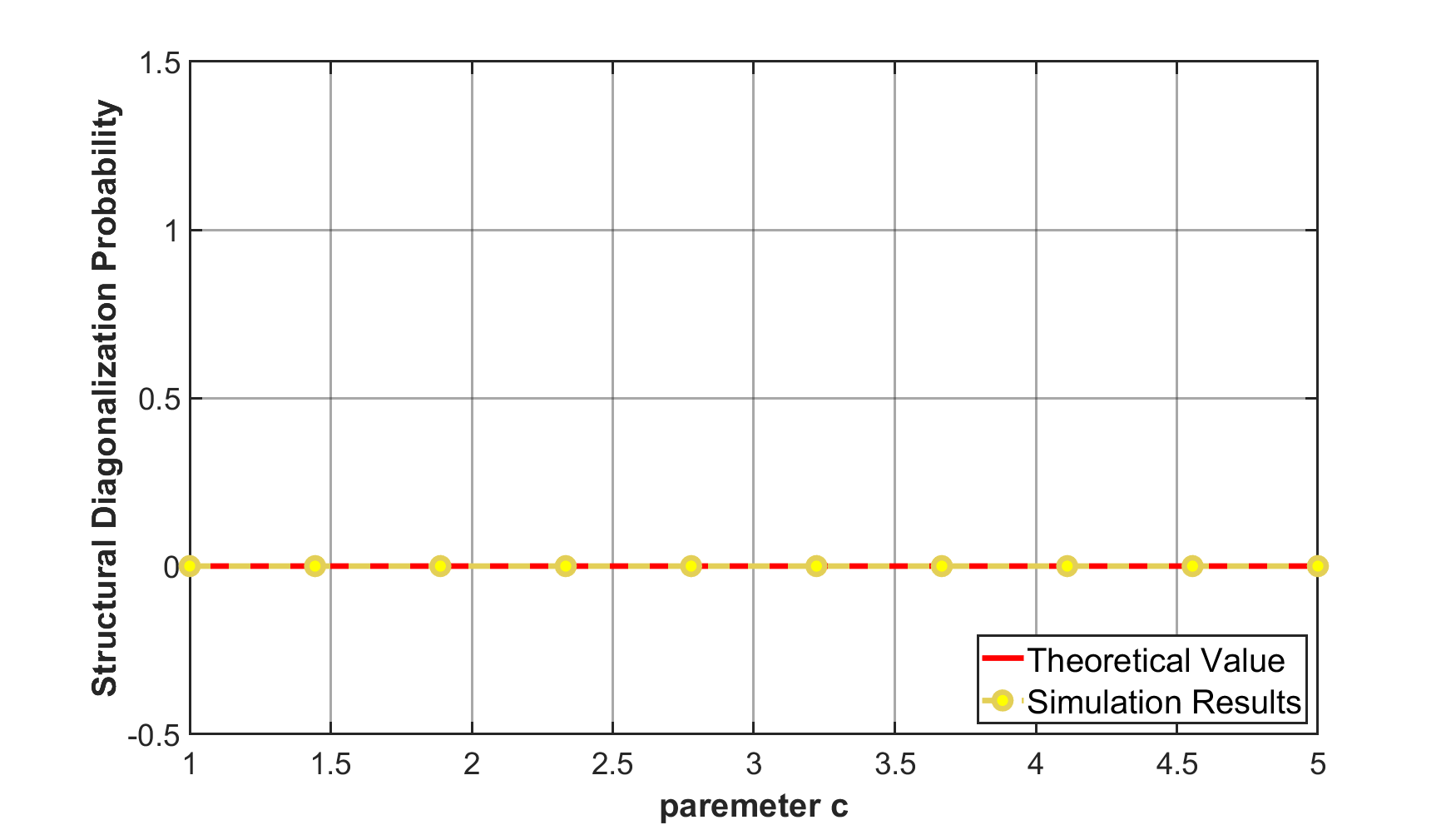}}
\caption{Comparison between simulation results and theoretical results for $p=\frac{c}{n}$ with fixed $n=4000$ and varying $c$}
\vspace{-1em}
\label{fig:n_n} 
\end{figure}

From Fig.~\ref{fig:n_n}, the simulation results align with theoretical predictions. These observations demonstrate that empirical results converge to theoretical expectations with high probability.

Now consider the asymptotic behavior of the simulation results compared to the theoretical results, and examine whether, under different connection probabilities, the simulation results gradually fall within the theoretical upper and lower bounds as the number of vertices n in the random graph increases.

Case 3: $p = \frac{\log n + c+o(1)}{n}$
\begin{itemize}
    \item  Fixed $n \in\{500,1000,1500,\dots ,3500,4000\}$ with $800$ samples per configuration
    \item $c \in \{-0.5,0,0.5,1,1.5,2\}$
    \item Theoretical bounds: \begin{equation*}
    e^{-2e^{-c}}\big(1 + 2e^{-c}+e^{-2c}\big) \leq \lim_{n\to \infty}\mathrm{P}(G_{n,p}\in  \mathcal{D}^{n})\leq 1-e^{-2e^{-c}} e^{-2c}.
\end{equation*}
\end{itemize}
\begin{figure}[htbp]
\centering
{\includegraphics[width=1\textwidth]{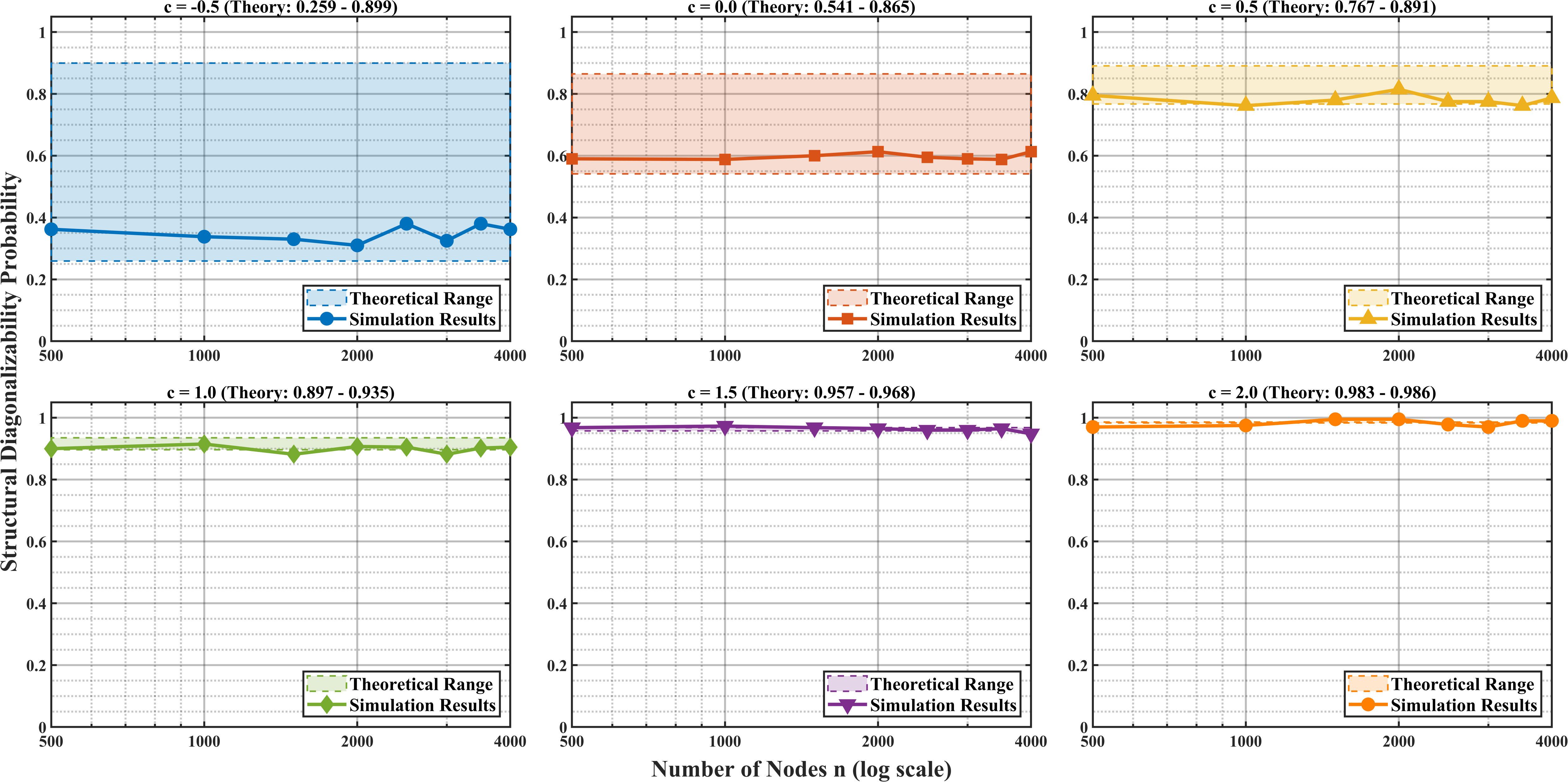}} 
\caption{Comparison between simulation results and theoretical upper/lower bounds for $p=\frac{\log n+c+o(1)}{n}$ with varying $c$ and $n$}
\vspace{-1em}
\label{fig:dc} 
\end{figure}

As shown in the six subfigures of Fig.~\ref{fig:dc}, regardless of the connection probability, as the number of vertices increases, the simulation results gradually stabilize around a certain value within the region between the theoretical upper and lower bounds. The fluctuations in the simulation curves diminish progressively, indicating that the simulation results asymptotically approach the theoretical values as the number of vertices increases.

\subsection{The simulation results of the $\mathcal{G}(n,p,q)$ random graph model}
The simulation framework contains three key parameters:
\begin{itemize}
    \item $p$: The probability of a non-self-loop edge, determined by $n$ and $c$
    \item $q$: The probability of a self-loop edge
    \item $n$: Number of vertices in generated graphs
  \end{itemize}

This experiment investigates the consistency between simulation results and theoretical bounds across different $c$ values under fixed network sizes.

Case1 : $p = \frac{\log n+c+o(1)}{n}$
\begin{itemize}
    \item Fixed $n =4000$ and $c=0$ with $800$ samples per configuration
    \item $q \in [0,1]$ sampled at 0.1 intervals for smooth curves
    \item Theoretical bounds:
    \begin{align*}
    e^{-2(1-q)e^{-c}}\big(1 + &2(1-q)e^{-c}+(1-q)^{2}e^{-2c}\big) \leq \lim_{n\to \infty}\mathrm{P}(G_{n,p}^q\in  \mathcal{D}^{n})\\&\leq
    1-(1-q)^{2}e^{-2(1-q)e^{-c}} e^{-2c}.
\end{align*}
\end{itemize}

\begin{figure}[htbp] 
\centering
\subfloat
  {\includegraphics[width=0.6\textwidth]{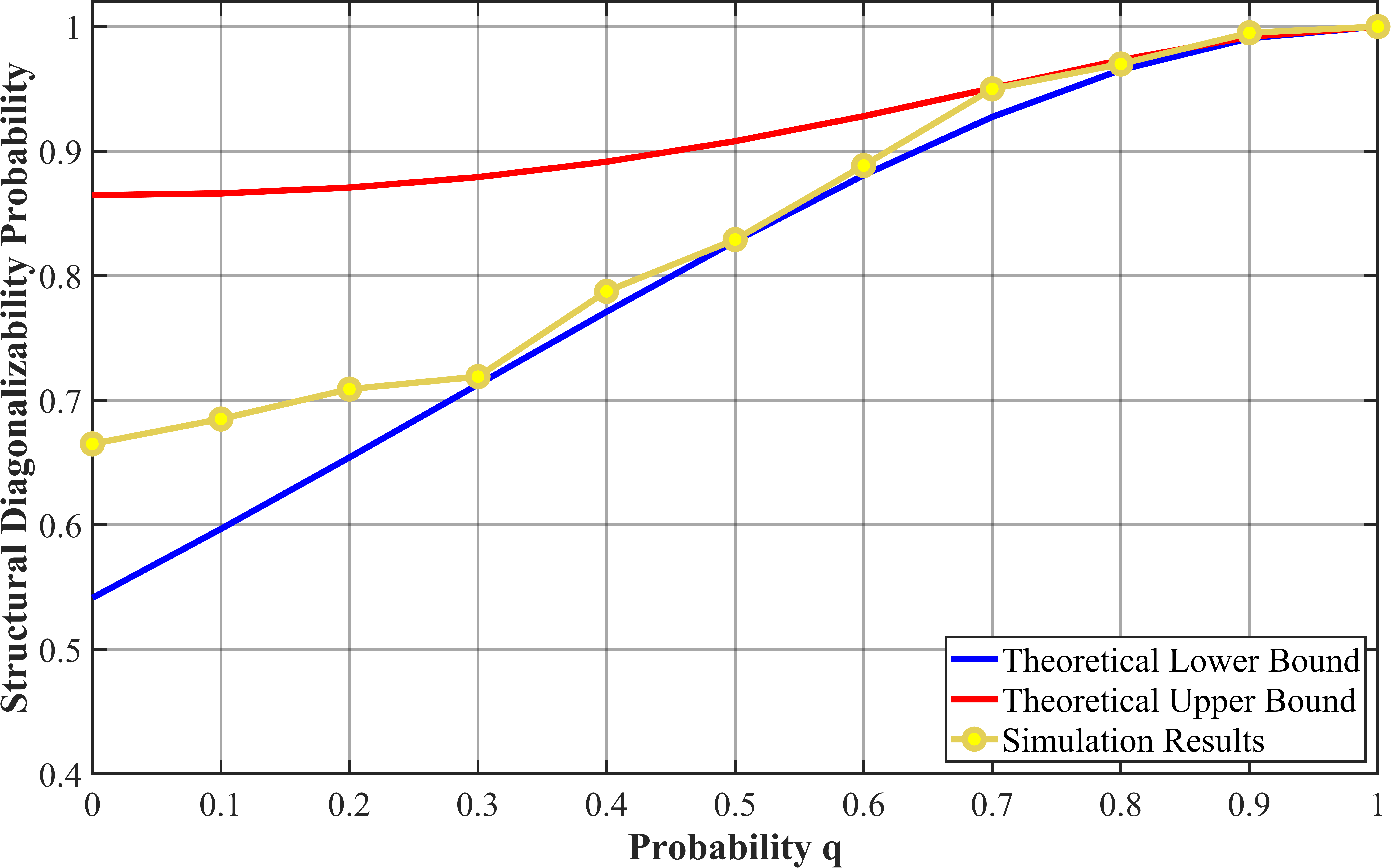}}
\caption{Comparison between simulation results and theoretical upper/lower bounds for $p=\frac{\log n+c+o(1)}{n}$ with fixed $n=4000$, $c=0$ and varying $q \in [0,1]$}
\vspace{-1em}
\label{fig:log_n_q_n} 
\end{figure}

Fig.~\ref{fig:log_n_q_n} demonstrates that when $c=0$, the structural diagonalizability probability for random digraphs ($n=4000$, $800$ instances per parameter set) is bounded from above and below by our theoretical predictions for all tested values of the self-loop probability $q$.

Next, we examine the empirical convergence of the results within the derived theoretical bounds as a function of increasing graph size and for different connection probabilities. The simulation settings are as follows. 

Case2 : $p = \frac{\log n+c+o(1)}{n}$
\begin{itemize}
    \item Fixed $n \in\{500,1000,1500,\dots ,3500,4000\}$ and $c=0$ with $800$ samples per configuration
    \item $q \in \{0,0.2,0.4,0.6,0.8,1\}$
    \item Theoretical bounds:
     \begin{align*}
    e^{-2(1-q)e^{-c}}\big(1 + &2(1-q)e^{-c}+(1-q)^{2}e^{-2c}\big) \leq \lim_{n\to \infty}\mathrm{P}(G_{n,p}^q\in  \mathcal{D}^{n})\\&\leq
    1-(1-q)^{2}e^{-2(1-q)e^{-c}} e^{-2c}.
\end{align*}
\end{itemize}

\begin{figure}[htbp]
\centering
  {\includegraphics[width=1\textwidth]{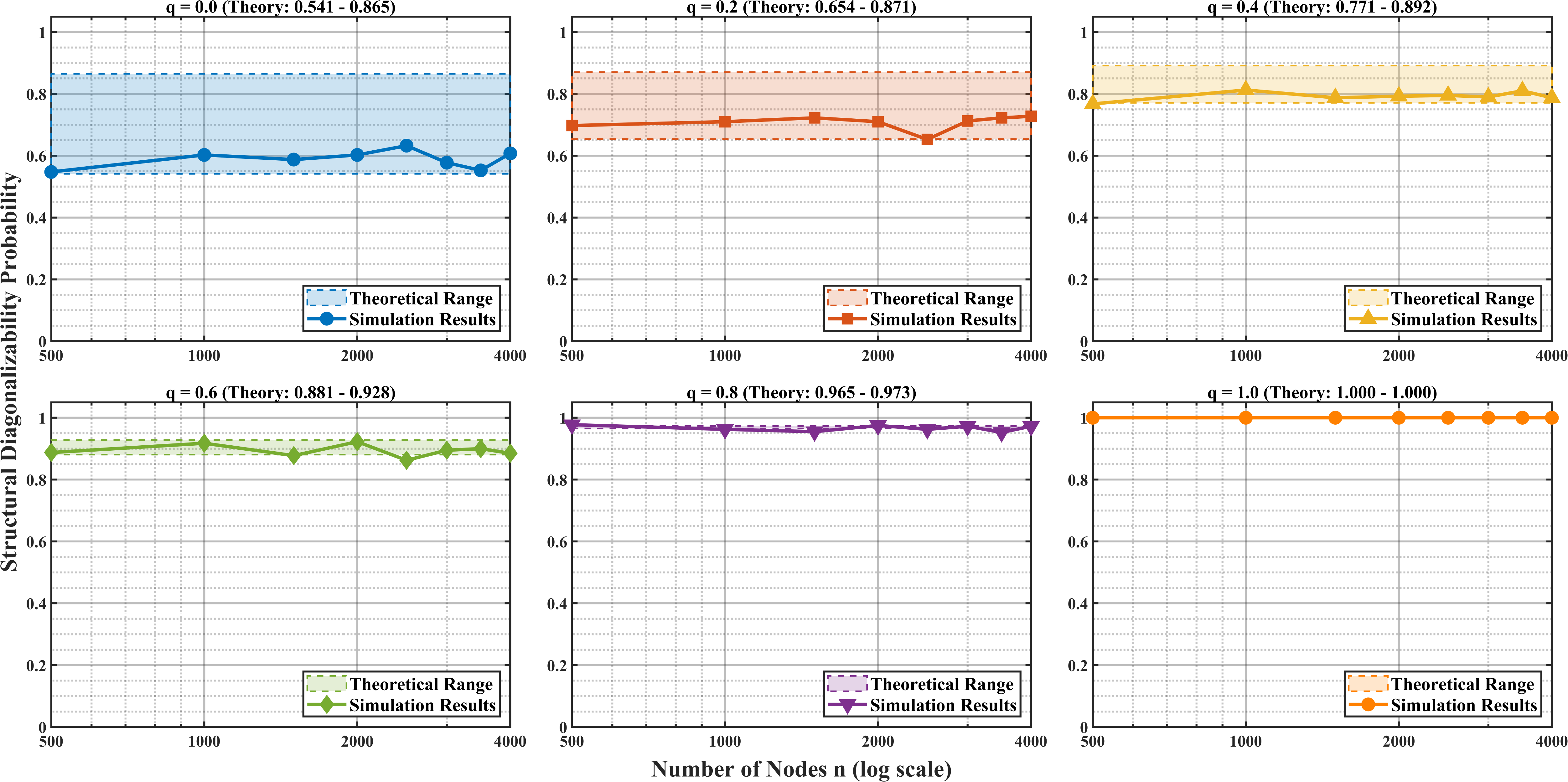}}
\caption{Comparison between simulation results and theoretical upper/lower bounds for $p=\frac{\log n+c+o(1)}{n}$ with fixed $c=0$ and varying $q\in \{0,0.2,0.4,0.6,0.8,1\}$}
\vspace{-1em}
\label{fig:dqc} 
\end{figure}

The six subfigures presented in Fig.~\ref{fig:dqc} demonstrate that under varying specifications of the non-self-loop connection probability, the simulation outputs consistently converge to values bounded by the theoretical limits as the network size increases. A marked attenuation in the oscillation amplitude of the simulation trajectories is observed, substantiating the asymptotic alignment between empirical simulations and theoretical predictions with increasing network size. This convergence provides empirical validation for the theoretical bounds governing the system's behavior.

The above simulation experiments have validated that for both the $\mathcal{G}(n,p)$ model and the $\mathcal{G}(n,p,q)$ random graph model, the asymptotic probabilities of structural diagonalizability derived through graph theory and probability theory are in complete agreement with the results of large-scale simulation experiments.

\section{Implications for structured systems} \label{sec-discussions}
We now briefly introduce structured systems and discuss implications of our results to the related problems within this framework. Consider a structured matrix quadruple $(\bar A,\bar B,\bar C, \bar  F)$, where $\bar A\in \{0,*\}^{n\times n}$, $\bar B\in \{0,*\}^{n\times m}$,  $\bar C \in \{0,*\}^{p\times n}$, and $\bar F \in \{0,*\}^{r\times n}$. Let $(A,B,C,F)$ be a realization of $(\bar A,\bar B,\bar C,\bar F)$.  A linear time-invariant system with parameters $(A,B,C,F)$ is described by the state equation and the output equation:
\begin{equation}\label{sys1}
    \begin{array}{c}
    \dot{x}(t)=Ax(t)+Bu(t)\\
    y(t)=Cx(t)\\
    z(t)=Fx(t)
    \end{array}
\end{equation}
where \(x(t)\in\mathbb{R}^{n}\) is the state vector, \(u(t)\in\mathbb{R}^{m}\) is the input vector, and \(y(t)\in\mathbb{R}^{p}\) is the output vector, $z(t)\in {\mathbb R}^r$ is a linear function of state. The quadruple $(\bar A,\bar B,\bar C, \bar F)$ encodes all such systems \eqref{sys1} with the same zero-nonzero patterns (termed structured systems).

System (\ref{sys1}) is said to be output controllable if for any $y_0,y_f\in {\mathbb R}^p$, there exists a finite $T>0$ and control input $u(t): [0,T)\to {\mathbb R}^m$ such that $y(0)=y_0$ and $y(T)=y_f$. System (\ref{sys1}) is said to be functionally observable if $z(0)$ can be uniquely determined from the input $u(t)|_{0\le t\le T}$ and output $y(t)|_{0\le t\le T}$ for some $T>0$. 
Recently, it has been found that output controllability (OC) and functional observability (FO) are dual to each other under certain conditions \cite{montanari2024popov}, which further enables the dual design of target controllers and functional observers \cite{montanari2024duality}. In general, numerical criteria for OC and FO involve high-order powers of the state matrix. When \(A\) is diagonalizable, however, simple Popov-Belevitch-Hautus tests exist \cite{zhang2023functional,montanari2024popov}. In this case, system \eqref{sys1} is OC if and only if ${\rm rank}(C[\lambda I-A,B])=p$ for all $\lambda\in \mathbb{C}$, and system \eqref{sys1} is FO if and only if ${\rm rank}([\lambda I-A^T,C^T,F^T]^T)={\rm rank}([\lambda I-A^T,C^T]^T)$ for all $\lambda\in \mathbb{C}$.  

A structured triple $(\bar A,\bar B,\bar C)$ is \emph{structurally output controllable (SOC)} if there exists a numerical realization that is output controllable. A triple $(\bar A,\bar C,\bar F)$ is \emph{structurally functionally observable (SFO)} if almost all its realizations $(A,C,F)$ are functionally observable~{\cite{montanari2022functional,zhang2023functional}}.

Verifying SOC for general structured systems is a longstanding open problem \cite{gao2014target,zhang2025generic}. The core criterion requires checking the generic rank of the output controllability matrix: $\operatorname{grank} (\bar{C}Q(\bar{A},\bar{B})) = p$, where $Q(\bar{A},\bar{B}) = [\bar{B}, \bar{A}\bar{B}, \dots, \bar{A}^{n-1}\bar{B}]$ involves high-order matrix powers. In contrast, for structurally diagonalizable systems, this verification is greatly simplified, admitting explicit graph-theoretic criteria involving dilations/contradictions that are verifiable in polynomial time {\cite{montanari2023target, zhang2025generic}}.
Similarly, verifying SFO in general structured systems requires checking the rank condition $$\operatorname{grank}\left[\begin{array}{c}O(\bar{A},\bar{C})\\
O(\bar{A},\bar{F})\end{array}\right] = \operatorname{grank} O(\bar{A},\bar{C}),$$where $O(\bar{A},\bar{C}) =$
$[\bar C^T,\bar A^T\bar C^T,...,(\bar C\bar A^{n-1})^T]^T$. Although checkable via weighted maximum matching algorithms with $O(n^3)$ complexity, the process is more efficient for structurally diagonalizable systems {\cite{montanari2022functional, zhang2025generic}}. For this class of systems, the condition simplifies to $\operatorname{grank} [\bar{A}^T, \bar{C}^T]^T = \operatorname{grank} [\bar{A}^T, \bar{C}^T, \bar{F}^T]^T$ combined with a graph connectivity condition, reducible to $O(n^2)$ complexity. 

Structural diagonalizability provides the structural counterpart of these numerical simplifications: if the interconnection graph is structurally diagonalizable, then almost all numerical realizations consistent with the graph admit a modal decomposition \cite{zhang2025generic}. This explains why structural diagonalizability can lead to simpler graph-theoretic criteria for SFO and SOC.

Early actuator-placement studies for SOC employed heuristic greedy strategies, such as the approach in
\cite{gao2014target}. More recent actuator and sensor placement methods admit polynomial-time solutions when SOC or SFO can be characterized by local graph properties, including dilations, contractions, and matching conditions \cite{zhang2025generic,zhang2025structural}. Structural diagonalizability provides a sufficient setting under which the original high-order criteria reduce to such local graph-theoretic and first-order generic-rank conditions {\cite{montanari2023target, zhang2025generic}}. Crucially, the corresponding actuator/sensor placement problems remain open in the general case \cite{gao2014target, zhang2025generic}.

Note that ${\rm P}(G_{n,0.5}\in \mathcal{D}^{n})$ is the proportion of structurally diagonalizable graphs among all digraphs on $n$ vertices. According to the results in Section 5.1.1, this probability ${\rm P}(G_{n,0.5}\in \mathcal{D}^{n})$ approaches $1$ as $n\to \infty$. This asymptotic behavior implies a significant practical insight: although SOC verification and the corresponding actuator placement problems are computationally intractable in the worst case, for almost all large-scale structured systems (i.e., with high probability under the random graph model), these problems can be solved efficiently in polynomial time. This bridges the gap between worst-case complexity and average-case one, highlighting the relevance of structural diagonalizability in practical system analysis. Similar results hold for the polynomial solvability of the sensor placement problem related to SFO. 

\section{Conclusions}
  To characterize the asymptotic probability of structural diagonalizability, we begin by examining its equivalent graph-theoretic conditions, which equate structural diagonalizability with the existence of a maximum matching that forms a disjoint cycle cover, along with isolation for all remaining vertices. To study the formation of disjoint cycles, we introduce the bipartite graph corresponding to the digraph and establish the relationship between disjoint cycle covers and the so-called consistent matchings. Furthermore, we propose some easily verifiable sufficient or necessary conditions for the structural diagonalizability of digraphs, laying a theoretical foundation for the subsequent characterization of upper and lower bounds on the asymptotic probability of structural diagonalizability in random graphs.
  Based on them, we provide a complete characterization of the asymptotic probability of structural diagonalizability for random digraphs with varying edge densities within the $\mathcal{G}(n,p)$ model. See \Cref{table:asymptotic_results} for details. 

As an extension of the classical $\mathcal{G}(n,p)$ model, we introduce a random digraph model $\mathcal{G}(n,p,q)$ with distinct probabilities for non-self-loop and self-loop edges. When the non-self-loop edge probability satisfies $p = \frac{\log n+c+o(1)}{n}$ and the self-loop probability is $q \in [0,1]$, we derive asymptotic upper and lower bounds for the probability of structural diagonalizability.

Our results also provide an affirmative answer to the conjecture made in
\cite{zhang2025generic}. These theoretical advances lead directly to practical applications: they significantly simplify the verification of system SFO and SOC, and they reduce the corresponding minimal controller and sensor design problems to computational tasks solvable in polynomial
time. Our probabilistic results may also provide insight into the local regularity of spectral-abscissa optimization for random structured closed-loop systems and controllability analysis of networked high-order systems with random topologies.

\appendix
\section{Proof of \cref{thm:main-result}}
    We first establish the upper bound of $\mathrm{P}(G_{n,p(n)}\in \mathcal{F}_k^n)$, where all edges (including the $n$ self-loop edges) appear independently with the same probability $p(n)$.
By \cref{def:Fk} and \cref{lem:The neighborhood set satisfies the condition.}, we obtain the following probability bound:
\begin{equation}
\mathrm{P}(G_{n,p(n)}\in \mathcal{F}_k^n) \leq 2 \cdot \binom{n}{k} \binom{n}{k-1} (1-p(n))^{(n-k+1)k} \cdot \left(\binom{k}{2} (p(n))^2\right)^{k-1}.
\end{equation}

The term $2$ represents that the vertex set $I$ satisfying the Condition~1 of $\mathcal{F}_k^n$ lies in the left or right part of $\mathcal{B}(G_{n,p(n)})$, and the term $\binom{n}{k}$ accounts for the selection of the $k$-vertex set $I$, and the term $\binom{n}{k-1}$ accounts for the choice of a $(k-1)$-vertex neighborhood set $N(I)$, and the term $(1-p(n))^{(n-k+1)k}$ accounts for the absence of edges between $I$ and $V \setminus N(I)$ and the term $\left(\binom{k}{2} (p(n))^2\right)^{k-1}$ accounts for the requirement that every vertex in $N(I)$ is adjacent to at least two vertices in $I$.

For integers $k \ge 2$, Stirling's approximation yields the two-sided inequality~\cite{belabbas2022stable},
\begin{equation}
\sqrt{2\pi}\,k^{k+0.5}e^{-k} \le k! \le e\,k^{k+0.5}e^{-k}.
\end{equation}
We proceed by expanding the binomial coefficients and applying these factorial estimates:
\begin{align}
\mathrm{P}(G_{n,p(n)}\in\mathcal{F}_k^n)
&\le 2 \cdot \frac{n!}{k!(n-k)!} \cdot \frac{n!}{(k-1)!(n-k+1)!} \nonumber \\
&\qquad \cdot \frac{k^{k-1}(k-1)^{k-1}}{2^{k-1}} p(n)^{2(k-1)} (1-p(n))^{(n-k+1)k} \nonumber \\
&\le 2 \cdot \frac{e n^{n+0.5}e^{-n}}{\sqrt{2\pi}k^{k+0.5}e^{-k}} \cdot \frac{e n^{n+0.5}e^{-n}}{\sqrt{2\pi}(k-1)^{k-0.5}e^{-k+1}} \nonumber \\
&\qquad \cdot \frac{1}{\sqrt{2\pi}(n-k)^{n-k+0.5}} \cdot \frac{1}{\sqrt{2\pi}(n-k+1)^{n-k+1.5}} \nonumber \\
&\qquad \cdot \frac{k^{k-1}(k-1)^{k-1}}{2^{k-1}} p(n)^{2(k-1)}(1-p(n))^{(n-k+1)k}.
\end{align}
Collecting the factorial terms gives
\begin{align*}
\mathrm{P}(G_{n,p(n)}\in\mathcal{F}_k^n)
&\le \frac{2 e^2 n^{2n+1}k^{k-1}(k-1)^{k-1}p(n)^{2(k-1)}(1-p(n))^{(n-k+1)k}}
{(2\pi)^2 k^{k+0.5}(n-k)^{n-k+0.5}(k-1)^{k-0.5}(n-k+1)^{n-k+1.5}2^{k-1}}.
\end{align*}

When $p(n)=\frac{\log n+c+o(1)}{n}$, we have
\begin{equation}
p(n)^{2(k-1)} = \frac{(\log n)^{2(k-1)}}{n^{2(k-1)}}\cdot
\frac{(\log n+c+o(1))^{2(k-1)}}{(\log n)^{2(k-1)}} .
\end{equation}
Define $C = \frac{\log n+c+o(1)}{\log n}$, which clearly satisfies $C=1$ for large $n$.  Also note that $\frac{2e^2}{(2\pi)^2}<1$.  Hence,
\begin{align}
&\mathrm{P}(G_{n,p(n)}\in \mathcal{F}_k^n) \nonumber \\
&\le \frac{n^{2n+1}(1-p(n))^{(n-k+1)k}}
{(k-1)^{2}(n-k)^{n-k+0.5}(n-k+1)^{n-k+1.5}2^{k-1}}
\frac{(\log n)^{2(k-1)}}{n^{2(k-1)}} \, C^{2(k-1)} \nonumber \\
&= \frac{n^{2n-2k+3}}
{(n-k)^{n-k+0.5}(n-k+1)^{n-k+1.5}}
\frac{(1-p(n))^{(n-k+1)k}(\log n)^{2(k-1)}}
{(k-1)^{2}2^{k-1}} \, C^{2(k-1)} \nonumber \\
&\le \frac{n^{2n-2k+3}}
{(n-k)^{n-k+0.5}(n-k)^{n-k+1.5}}
\frac{(1-p(n))^{(n-k+1)k}(\log n)^{2(k-1)}}
{(k-1)^{2}2^{k-1}} \, C^{2(k-1)} \nonumber \\
&= \left(\frac{n}{n-k}\right)^{\!2n-2k+2}\,
\frac{n\,(1-p(n))^{(n-k+1)k}(\log n)^{2(k-1)}}
{(k-1)^{2}2^{k-1}} \, C^{2(k-1)} \nonumber \\
&= \left(1+\frac{k}{n-k}\right)^{\!2n-2k+2}\,
\frac{n\,(1-p(n))^{(n-k+1)k}(\log n)^{2(k-1)}}
{(k-1)^{2}2^{k-1}} \, C^{2(k-1)} \nonumber \\
&\le \Biggl(\left(1+\frac{k}{n-k}\right)^{\!\frac{2n-2k+2}{k}}
n^{\frac{1}{k}}(1-p(n))^{n-k+1}(\log n)^{2} \, C^{2}\Biggr)^{\!k}.
\end{align}
For $n\ge 3$ and $2\le k\le n-1$, the inequality $\frac{2n-2k+2}{k}\le \frac{4(n-k)}{k}$ holds.  Consequently,
\begin{align*}
\mathrm{P}(G_{n,p(n)}\in\mathcal{F}_k^n)
&\le \Biggl(\left(1+\frac{k}{n-k}\right)^{\!\frac{4(n-k)}{k}}
n^{\frac{1}{k}}(1-p(n))^{n-k+1}(\log n)^{2} \, C^{2}\Biggr)^{\!k}.
\end{align*}
Using the standard bound $(1+\frac{1}{x})^{x}\le e$, we obtain
\begin{equation}
\mathrm{P}(G_{n,p(n)}\in\mathcal{F}_k^n) \le
\Bigl(e^{4}\, n^{\frac{1}{k}}(1-p(n))^{n-k+1}(\log n)^{2} \, C^{2}\Bigr)^{\!k}.
\end{equation}

A Taylor expansion yields
\begin{equation}
1-p(n)=e^{\log(1-p(n))}=e^{-p(n)} =
e^{-\frac{\log n+c+o(1)}{n}} =
n^{-\frac{1}{n}\bigl(1+\frac{c}{\log n}+o(1)\bigr)},
\end{equation}
Substituting this gives
\begin{align}
\mathrm{P}(G_{n,p(n)}\in\mathcal{F}_k^n)
&\le \Bigl(e^{4}(\log n)^{2} \, C^{2} \,
n^{-\frac{n-k+1}{n}\bigl(1+\frac{c}{\log n}+o(1)\bigr)+\frac{1}{k}}\Bigr)^{\!k} \nonumber \\
&= \Bigl(e^{4}(\log n)^{2} \, C^{2} \,
n^{-\bigl(\frac{n-k+1}{n}(1+\frac{c}{\log n}+o(1))-\frac{1}{k}\bigr)}\Bigr)^{\!k}.
\end{align}
Recalling that $C=\frac{\log n+c+o(1)}{\log n}$, we have
\begin{equation}
\mathrm{P}(G_{n,p(n)}\in\mathcal{F}_k^n) \le
\Bigl(e^{4}\bigl(\log n+c+o(1)\bigr)^{2} \,
n^{-\bigl(\frac{n-k+1}{n}(1+\frac{c}{\log n}+o(1))-\frac{1}{k}\bigr)}\Bigr)^{\!k}.
\end{equation}

For any $c\in\mathbb{R}$ or $c=\pm\omega(1)$, and for $n$ sufficiently large so that $p(n)\ge0$, we observe
\begin{equation}
\frac{n-k+1}{n}\Bigl(1+\frac{c}{\log n}+o(1)\Bigr)-\frac{1}{k}
\;\ge\; \frac{n-k+1}{n}-\frac{1}{k}
\;=\; \frac{(k-1)(n-k)}{nk},
\end{equation}
Since $2\le k\le n-1$, the right‑hand side is positive.  Therefore,
\begin{equation}
\mathrm{P}(G_{n,p(n)}\in\mathcal{F}_k^n) \le
\Biggl(e^{4}\,
\frac{\bigl(\log n+c+o(1)\bigr)^{2}}
{n^{\frac{n-k+1}{n}(1+\frac{c}{\log n}+o(1))-\frac{1}{k}}}\Biggr)^{\!k}.
\end{equation}

In conclusion, for $c\in\mathbb{R}$ or $c=\pm\omega(1)(c\geq -\log n)$ and $2\le k\le\lceil n/2\rceil$,
\begin{equation}
e^{4}\,\frac{\bigl(\log n+c+o(1)\bigr)^{2}}
{n^{\frac{n-k+1}{n}(1+\frac{c}{\log n}+o(1))-\frac{1}{k}}}=o(1).
\end{equation}
Summing over $k$ and using the geometric series formula, we obtain
\begin{equation}\label{equ:bound}
\lim_{n\to\infty}\sum_{k=2}^{\lceil n/2\rceil}\mathrm{P}(G_{n,p(n)}\in\mathcal{F}_k^n)=o(1).
\end{equation}

 The events \(\mathcal F_k^n\), \(1\le k\le \lceil n/2\rceil\), describe the possible Hall-type obstruction events contributing to \(\Gamma^n\). Eq.~\eqref{equ:bound} shows that, under the model \(\mathcal G(n,p)\), the total contribution of all obstruction events with \(k\ge2\) is asymptotically negligible.

We next extend the above estimate to the model
\(\mathcal G(n,p,q)\). To this end, we compare the probabilities of the events \(\mathcal F_k^n\), \(k\geq2\), with those in the corresponding model without self-loops through a coupling argument.

Construct $G_{n,p(n)}^{q(n)}$ and a reference graph $G^*$ on the same probability space as follows:
\begin{enumerate}
    \item All non-loop edges are the same in both graphs (same Bernoulli $p(n)$ realizations),
    \item The self-loop edges in $G_{n,p(n)}^{q(n)}$ are independent Bernoulli trials with probability $q(n)$,
    \item The self-loop edges in $G^*$ are forced to be absent (i.e., probability $0$).
\end{enumerate}

Then almost surely
\[
E(G^*) \subseteq E(G_{n,p(n)}^{q(n)}).
\]

Since adding edges can only increase (or keep unchanged) the size of any neighborhood set $|N(I)|$,
the event $\mathcal{F}_k^n$ (for $k \ge 2$) becomes strictly harder to satisfy when more edges are present.
In particular,
\[
\mathrm{P}\bigl( G_{n,p(n)}^{q(n)} \in \mathcal{F}_k^n \bigr) \;\le\; \mathrm{P}\bigl( G^* \in \mathcal{F}_k^n \bigr).
\]

Observe that the difference between $G^*$ and $G_{n,p(n)}$ lies solely in the absence of self-loop edges in $G^*$. This affects at most $k$ possible edges incident to the set $I$ (namely, the self-loop edges from vertices in $I$ to their counterparts in the opposite part), since $|I| = k$ is fixed or grows slowly compared to $n$.
Thus, when $p(n)=\frac{\log n+c+o(1)}{n}$, we have
\[
\sum_{k=2}^{\lceil n/2 \rceil} \mathrm{P}\bigl( G^* \in \mathcal{F}_k^n \bigr)\leq \sum_{k=2}^{\lceil n/2 \rceil} 2  \binom{n}{k} \binom{n}{k-1} (1-p(n))^{(n-k+1)k}  \left(\binom{k}{2} (p(n))^2\right)^{k-1}=o(1),
\]which we shall explain as follows. The term $(1-p(n))^{(n-k+1)k}$ in the union-bound expression represents the probability of no edges from $I$ to $V \setminus N(I)$. In $G^*$, the true probability of no such edges is at least as large as (and typically larger than) in $G_{n,p(n)}$, because some potential self-loop edges are forced absent. The ratio between the two probabilities is bounded by
    \[
    \left( \frac{1-0}{1-p(n)} \right)^{O(k)} = \exp\bigl( O(k p(n)) \bigr) = \exp\bigl( O(k \cdot \frac{\log n}{n}) \bigr) = 1 + O\left( \frac{k \log n}{n} \right) = 1 + o(1),
    \]
    since $k$ is at most $\lceil n/2 \rceil$ but the exponential decay in the final bound is much stronger than any polynomial factor in $n$. Thus, replacing the term with its value from $G_{n,p(n)}$ yields a conservative upper bound.
The term $\left( \binom{k}{2} p(n)^2 \right)^{k-1}$ provides a crude upper bound on the probability that every vertex in $N(I)$ is adjacent to at least two vertices in $I$.  In $G^*$ (where self-loop edges are absent), the actual probability that every vertex in $N(I)$ has at least two neighbors in $I$ can only decrease or remain unchanged compared to $G_{n,p(n)}$. Therefore, the crude upper bound term 
$\left( \binom{k}{2} p(n)^2 \right)^{k-1}$ remains a valid (though looser) upper bound on the actual probability under $G^*$.

Therefore, when $p(n)=\frac{\log n+c+o(1)}{n}$, we have 
\begin{align*}
    &\sum_{k=2}^{\lceil n/2 \rceil} \mathrm{P}\bigl( G_{n,p(n)}^{q(n)} \in \mathcal{F}_k^n \bigr)\leq \sum_{k=2}^{\lceil n/2 \rceil} \mathrm{P}\bigl( G^* \in \mathcal{F}_k^n \bigr)=o(1).
\end{align*}
That is, under the same $p(n) = \frac{\log n + c + o(1)}{n}$, but with self-loop probability $q(n) \in [0,1]$,
we still have
\begin{equation*}
\lim_{n\to\infty}\mathrm{P}(G_{n,p(n)}^{q(n)}\in\Gamma^n)=\sum_{k=1}^{\lceil n/2 \rceil} \mathrm{P}\bigl( G_{n,p(n)}^{q(n)} \in \mathcal{F}_k^n \bigr)
= \lim_{n\to\infty}\mathrm{P}(G_{n,p(n)}^{q(n)}\in\mathcal{F}_1^n)+o(1).
\end{equation*}\hfill \proofbox

\section{Proof of \cref{thm:SDB=0}}

Let $I_{u,v}$ be the indicator that $(u\to v)$ is an isolated edge, and let
\[
N=\sum_{\substack{u\neq v\\ u,v\in\{1,\cdots, n\}}} I_{u,v}
\]
be the total number of isolated edges.

First compute the expectation. For a fixed ordered pair $(u,v)$ with $u \neq v$, the edge $u \to v$ is isolated if it exists and all other edges that would violate the degree conditions are absent. Specifically, the following conditions must hold: the edge $u \to v$ exists (probability $p(n) = \frac{c}{n}$); there are no loops at $u$ or $v$ (each with probability $1-p(n)$); there is no reverse edge $v \to u$ (probability $1-p(n)$); and for every other vertex $w \notin \{u,v\}$, there are no edges $u \to w$, $w \to u$, $v \to w$, or $w \to v$ (each with probability $1-p(n)$, and there are $n-2$ such vertices for each type). Since all these edges are distinct, their presence or absence are independent events. Multiplying the probabilities yields $\mathbb{E}I_{u,v} = p(n)(1-p(n))^{4n-5}$. Since $p(n)=\frac{c}{n}$ and $(1-p(n))^{4n-5}=e^{-(4n-5)p(n)+O(np(n)^2)}=e^{-4c+O(1/n)}$. Hence,
\[
\mathbb{E}N=\sum_{u\neq v}\mathbb{E}I_{u,v}=n(n-1)\cdot\frac{c}{n}e^{-4c}(1+o(1))
      =c\,e^{-4c}\,n\,(1+o(1))\longrightarrow\infty\quad (n\to\infty).
\]

Now we bound the variance. Write
\[
\operatorname{Var}(N)=\sum_{u\neq v}\operatorname{Var}(I_{u,v})
                    +\sum_{(u,v)\neq(u',v')}\operatorname{Cov}(I_{u,v},I_{u',v'}),
\]where ${\rm Var}(I_{u,v})$ takes the variance and ${\rm Cov}(I_{u,v},I_{u',v'})$ the covariance.
For the first sum, since $I_{u,v}$ is a Bernoulli variable,
\[
\operatorname{Var}(I_{u,v})=\mathbb{E}(I_{u,v})^2-(\mathbb{E}I_{u,v})^{2}=\mathbb{E}I_{u,v}-(\mathbb{E}I_{u,v})^{2}
                        \le\mathbb{E}I_{u,v},
\]
and therefore
\[
\sum_{u\neq v}\operatorname{Var}(I_{u,v})\le\sum_{u\neq v}\mathbb{E}I_{u,v}
                                        =\mathbb{E}N.
\]

For the covariance terms, we distinguish two cases.
\begin{enumerate}
\item If  $\{u,v\}$ and $\{u',v'\}$ are disjoint, then the indicators $I_{u,v}$ and $I_{u',v'}$ depend on disjoint sets of possible edges and are independent; hence $\operatorname{Cov}(I_{u,v},I_{u',v'})=0$.

\item If $\{u,v\}\cap\{u',v'\}\neq\emptyset$, then the two events $\{I_{u,v}=1\}$ and $\{I_{u',v'}=1\}$ cannot occur simultaneously. Indeed, a vertex that belongs to both ordered pairs would have to satisfy contradictory degree conditions required by an isolated edge. Consequently $I_{u,v}I_{u',v'}=0$, and
\[
\operatorname{Cov}(I_{u,v},I_{u',v'})
   =\mathbb{E}I_{u,v}I_{u',v'}-\mathbb{E}I_{u,v}\mathbb{E}I_{u',v'}
   =-\mathbb{E}I_{u,v}\mathbb{E}I_{u',v'}\le0.
\]
\end{enumerate}
Thus, every covariance is non‑positive, whence
\[
\sum_{(u,v)\neq(u',v')}\operatorname{Cov}(I_{u,v},I_{u',v'})\le0.
\]

Combining the estimates, we obtain
\[
\operatorname{Var}(N)\le\mathbb{E}N,
\qquad\text{so } \operatorname{Var}(N)=O\bigl(\mathbb{E}N\bigr).
\]

Finally apply Chebyshev's inequality \cite[Lemma 22.3]{frieze2015introduction} and use $\mathbb{E}N\to\infty$ together with $\operatorname{Var}(N)=O(\mathbb{E}N)$,
\[
\frac{\operatorname{Var}(N)}{(\mathbb{E}N)^{2}}
   =O\!\left(\frac{1}{\mathbb{E}N}\right)\longrightarrow0.
\]
Hence
\[
{\rm P}(N=0)\le{\rm P}\bigl(|N-\mathbb{E}N|\ge\mathbb{E}N\bigr)\le\frac{\operatorname{Var}(N)}{(\mathbb{E}N)^{2}}
               \longrightarrow0.
\]

By \cref{lem:gulibian}, a structurally diagonalizable graph cannot contain isolated edges. We therefore obtain
\begin{equation*}
\lim_{n \to \infty} \mathrm{P}(G_{n,p(n)}\in\mathcal{D}^{n})\leq \lim_{n\to\infty}{\rm P}(N=0)=0.
\end{equation*}\hfill \proofbox

\bibliographystyle{siamplain}
\bibliography{references}

\end{document}